\documentclass[11pt,letterpaper]{article}
\usepackage[letterpaper, margin=1in]{geometry}
\usepackage{graphicx,color}
\usepackage[cp850]{inputenc}
\usepackage[T1]{fontenc}
\usepackage[english]{babel}
\usepackage{rotating}
\usepackage[f]{esvect}
\usepackage{amsmath, amssymb}
\usepackage{amsthm}
\usepackage{rotating}
\usepackage{lscape}
\usepackage{mathrsfs}
\usepackage{mathtools}
\usepackage[toc,page]{appendix}
\usepackage{hyphenat}
\usepackage{hyperref}
\usepackage{authblk}

\usepackage[normalem]{ulem}

\def\tablenotes{\bgroup\parfillskip=0pt plus 1fil
\leftskip=0pt\relax \rightskip=0pt
\vskip2pt\footnotesize}
\def\endtablenotes{\vskip1pt\egroup}

\newtheorem{theorem}{Theorem}[section]

\newtheorem{corollary}[theorem]{Corollary}
\newtheorem{lemma}[theorem]{Lemma}
\newtheorem{remark}[theorem]{Remark}
\newtheorem{definition}[theorem]{Definition}

\newcommand{\TimeDeriv}{\frac{\textrm{d}}{\textrm{dt}}}

\renewcommand{\leq}{\leqslant}
\renewcommand{\geq}{\geqslant}
\renewcommand{\d}{\mathrm{d}}
\DeclarePairedDelimiter\ceil{\lceil}{\rceil}
\DeclarePairedDelimiter\floor{\lfloor}{\rfloor}

\newcommand{\Poisson}{\mathrm{Poisson}}
\newcommand{\Erlang}{\mathrm{Erlang}}
\newcommand{\GammaDist}{\mathrm{Gamma}}
\newcommand{\Var}{\mathrm{Var}}

\usepackage{amsopn}

\renewcommand{\phi}{\varphi}
\renewcommand{\epsilon}{\varepsilon}
 
\allowdisplaybreaks

 \usepackage[square,authoryear,sort]{natbib}

\usepackage{graphicx}
\usepackage{overpic}
\usepackage{tikz}
\numberwithin{equation}{section}

\makeatletter
\renewcommand{\@biblabel}[1]{#1\hfill \hspace{-0.1cm}}
\makeatother

  \graphicspath{ {./Figures/} }
 
\begin{document}
 \title{Numerical methods and hypoexponential approximations for gamma distributed delay differential equations }
 
 \author[1,*]{Tyler Cassidy}
 \author[2]{Peter Gillich}
 \author[2,3]{Antony R. Humphries}
 \author[1]{Christiaan H. van Dorp} 
 
 \affil[1]{Theoretical Biology and Biophysics (T-6), Los Alamos National Laboratory, Los Alamos NM, USA}
 \affil[2]{Department of Mathematics and Statistics, McGill University, Montr\'{e}al, Quebec, Canada}
 \affil[3]{Department of Physiology,  McGill University, Montr\'{e}al, Quebec, Canada } 
 \affil[*]{Corresponding author: tcassidy@lanl.gov}
 
\maketitle

\begin{abstract}
Gamma distributed delay differential equations (DDEs) arise naturally in many modelling applications. However, appropriate numerical methods for generic Gamma distributed DDEs are not currently available. Accordingly, modellers often resort to approximating the gamma distribution with  an Erlang distribution and using the linear chain technique to derive an equivalent system of ordinary differential equations. In this work, we develop a functionally continuous Runge-Kutta method to numerically integrate the gamma distributed DDE and perform numerical tests to confirm the accuracy of the numerical method. As the functionally continuous Runge-Kutta method is not available in most scientific software packages, we then derive hypoexponential approximations of the gamma distributed DDE. Using our numerical method, we show that while using the common Erlang approximation can produce solutions that are qualitatively different from the underlying gamma distributed DDE, our hypoexponential approximations do not have this limitation. Finally, we implement our hypoexponential approximations to perform statistical inference on synthetic epidemiological data. 

\end{abstract}
% \newpage

\section{Introduction}

Gamma distributed delay differential equations (DDEs), generically of the form
\begin{equation}
\left.
\begin{aligned}
\TimeDeriv X(t) & = F \left(X(t), \int_0^{\infty} X(t-s)g_a^j(s) \d s \right)\\
X(s) & = \psi (s) \quad s < t_0.
\end{aligned}
\right \}
\label{Eq:DDEIVP}
\end{equation}
have been extensively used in mathematical biology, epidemiology and pharmacometric modelling \citep{Smith2011,Cassidy2020a,DeSouza2017,Hurtado2019,Champredon2018,Hu2018,Ando2020}. These models describe the influence of the past on the current state through the convolution integral 
\begin{align*}
(X*g_a^j)(t) = \int_0^{\infty} X(t-s) g_a^j(s) \d s 
\end{align*}
where $g_a^j(s)$ is the probability density function (PDF) of the gamma distribution. The initial value problem~\eqref{Eq:DDEIVP} is equipped with initial data in the form of the history function $\psi$. Typically, $\psi \in L_1(\mathbb{P})$ where $\mathbb{P}$ is a probability measure \citep{Hale1993}. The Radon-Nikodym derivative of $\mathbb{P}$ with respect to Lebesgue measure is the PDF $g_a^j(s)$ given by
\begin{align*}
g_a^j (s) = \frac{a^j}{\Gamma(j)}s^{j-1} \exp(-a s).
\end{align*}
which is parameterized using the shape and scale parameters, $j$ and $a$, respectively. While both these parameters can be positive reals, in most applications, $j$ is restricted to integer values and \eqref{Eq:DDEIVP} is thus an Erlang distributed DDE. These Erlang distributed DDEs can then be reduced to an equivalent system of ordinary differential equations (ODEs) through the linear chain technique \citep{DeSouza2017,Vogel1961}. A major impediment to the implementation of the more general gamma distributed DDE is the lack of appropriate numerical techniques for their simulation \citep{Diekmann2020c,Diekmann2017,Breda2016}. Here, we address this by developing a functionally continuous Runge-Kutta (FCRK) method to simulate \eqref{Eq:DDEIVP} and deriving a finite dimensional approximation that is more accurate than the common Erlang approximation.

Currently existing numerical tools for the simulation and study of DDEs have only recently started to  include equations with infinite delay \citep{Diekmann2020c,Diekmann2017,Cassidy2018a}. In particular, work towards these goals includes pseudo-spectral techniques \citep{Gyllenberg2018,Diekmann2020c}, and the development of ODE approximations of the gamma distributed DDE without enforcing $j \in \mathbb{N}$ \citep{Krzyzanski2019,Koch2015}.  \citet{Krzyzanski2019} used the binomial theorem to develop an ODE approximation of a generic gamma distributed DDE. However, this approximation relies on truncating the infinite series expansion of the probability density function of the gamma distribution at some finite value. While \citet{Krzyzanski2019} does derive explicit error bounds dependent on the number of terms $M$ in the series expansion, the artificial truncation of the convolution integral ensures that the numerical approximation is not consistent. In a related work focused on lifespan distributions, \citet{Koch2015} impose a fixed upper bound for the lifespan duration, then subdivide the interval of possible lifespan durations into $m$ sub-compartments. The populations in each sub-compartment are weighted according to the probability of a lifespan of that length when calculating the total population size. Once again, this method requires the modeller to determine a fixed upper bound of the lifespan duration, and does not capture the full dynamics of the infinite delay DDE. The FCRK method developed in this work explicitly computes the improper convolution integral and eliminates the requirement that modellers impose artificial upper bounds. 
 We make a number of modifications to existing numerical methods for DDEs to derive our FCRK method for the infinite delay problem \eqref{Eq:DDEIVP}. The main difficulty in adapting existing FCRK methods to infinite delay problems is the accurate evaluation of the convolution integral. To address this, we introduce a change of variable that maps the semi-infinite domain of integration to a compact interval. It is then possible to utilize existing Newton-Coates methods to numerically calculate the transformed convolution integral and thus evaluate the right-hand side of \eqref{Eq:DDEIVP}. To ensure the accuracy of the FCRK method, we derive order conditions on the quadrature method, and demonstrate the accuracy of our FCRK method through a number of examples. 

Inspired by the lack of existing appropriate numerical methods for problems such as \eqref{Eq:DDEIVP}, there has also been considerable interest in approximating infinite delay DDEs by forms that are more convenient for simulation \citep{Cassidy2018a,Hurtado2019,Diekmann2017,Diekmann2020,Krzyzanski2019,Koch2015}. The most well-known of these is the previously mentioned linear chain technique, wherein modellers often make the simplifying assumption that $j \in \mathbb{N}$ when implementing gamma distributed DDE models. However, this assumption imposes constraints on the sample mean and variance of the delayed process. Typically, for a general gamma distributed random variable with mean $\tau$ and variance $\sigma^2$, modellers impose $ j = [\tau^2/\sigma^2]$ where $[x]$ rounds $x$ to the nearest integer with $[0.5]=1$ \citep{Cassidy2019,Jenner2020}. As a result, it is only possible to fit one of these statistics with an Erlang distribution, as the system
\begin{align*}
   \tau = \frac{[\tau^2/\sigma^2]}{a}  \quad \textrm{and} \quad \sigma^2 = \frac{[\tau^2/\sigma^2]}{a^2}
\end{align*}
only admits a solution if $[\tau^2/\sigma^2] = \tau^2/ \sigma^2$ which corresponds to $j \in \mathbb{N}.$

In light of this limitation of the Erlang approximation, recent work has explored using phase type distributions to approximate generic distributed DDEs. These phase type distributed DDEs are then reduced to a system of ODEs through a variant of the linear chain technique \citep{Hurtado2020,Hurtado2019}. These phase type distributions approximate the underlying distribution either by minimizing some distance measurement between distributions or by matching moments of the underlying distribution. In this work, we take the latter approach when developing a novel hypoexponential approximation of the generic gamma distributed DDE~\eqref{Eq:DDEIVP}. Existing work has identified the ``reachable'' bounds in moment space and the minimal number of phases required to match the first three moments of the underlying distribution \citep{Bobbio2005,Johnson1990,Johnson1989}. However, in general, it may not be possible to match the first three moments, and even if it is possible, the required number of phases can be arbitrarily high \citep{Osogami2006}. To avoid the aforementioned complications, we only match the first two moments of the underlying gamma distribution. We match the first two moments without imposing any restrictions on their values, and the parameters of our approximating phase type distribution are entirely determined by the mean and variance of the underlying distribution. Moreover, our approximation is exact if the underlying distribution is Erlang. 

To achieve this, we derive explicit rates of the hypoexponenetial approximation to match the first two moments of a given gamma distribution, and then derive the equivalent system of ODEs. These ODE models are simple to study numerically and have the added benefit of being easy to implement and explain for non-mathematical users. Further, we leverage our FCRK method to simulate \eqref{Eq:DDEIVP} and thus explicitly evaluate the accuracy of our hypoexponenetial approximation, which has not been done in prior work. As we will show, our approximation more accurately captures the dynamics of the underlying gamma distributed DDE than the Erlang approximation obtained by setting $ j = [\tau^2/\sigma^2]$. Finally, we utilize the explicit formulas for the rates in our hypoexponenetial approximation to show that it is exact if the gamma distribution has an integer shape parameter $j$.  

Finally, we apply our hypoexponenetial approximation to the problem of statistical inference. Erlang delay models are often used in epidemiology of infectious diseases \citep{Rozhnova2021,Champredon2018,Sanche2020,Greenhalgh2019}, but also in many other fields \citep{Smith2011,Cassidy2018,DeSouza2017,Cassidy2020a}. 
A common approach is to define an epidemiological model in terms of ODEs, and estimate parameters by fitting the model to disease incidence data. One key parameter, the basic reproduction number $R_0$, is closely related to the generation interval (or it's proxy the serial interval) and the initial exponential growth rate of the epidemic, and the incubation period $T_E$ and time to recovery $T_I$ \citep{roberts_model_2007}. This relation depends not only on the mean infectious period $\mathbb{E}[T_I]$ and incubation period $\mathbb{E}[T_E]$, but also on the distribution of these periods. Therefore, making invalid assumptions about this distribution can lead to spurious estimates of e.g.\ $R_0$.

In practice, the random times $T_E$ and $T_I$ are often assumed to be $\Erlang$-distributed so that the DDEs can simply be implemented as ODEs using the linear chain technique. This is convenient because ODE models are easy to implement in commonly used software packages for statistical inference \citep{Carpenter2017}, whereas support for DDE models is much less common, and often restricted to an Erlang distributed or fixed delay \citep{monolix_2019,Raue2009,Raue2014}. Apart from convenience, there is no reason to assume that the distributions of $T_E$ and $T_I$ should be $\Erlang$ instead of gamma distributions. The hypoexponential approximation of the gamma distribution proposed here still allows for an ODE approximation of the DDE model, but removes the need to assume that the shape parameter is an integer.

Another purely practical reason for using the hypoexponential approximation of the gamma distribution instead of an $\Erlang$ distribution is that estimating the integer shape parameter $j$ of the Erlang distribution can be inconvenient in some software packages. For instance, in the commonly used package for Bayesian inference Stan \citep{Carpenter2017}, the estimated parameters have to be real-valued due to the limitations imposed by the Hamiltonian Monte-Carlo method. This means that in order to estimate an integer valued $j$, one would have to repeat the analysis for multiple values of $j$, and compare models with, ideally Bayes factors, or more practically, with information criteria as 
LOO-IC or WAIC \citep{Vehtari2017} to combine the posterior samples from these separate model fits. This extra step and the resulting extra computation can be avoided if $j$ is allowed to be real-valued, as in the approximation derived in this work. Consequently, we implement the hypoexponential approximation in Stan and use the resulting ODE system for statistical inference of epidemiological parameters and evidence synthesis.

The remainder of the article is structured as follows. We begin by developing a numerical method to simulate the general gamma distributed DDE by using the theory of functionally continuous Runge-Kutta methods to address the overlapping in the convolution integral in Section~\ref{Sec:FCRK}. Next, we develop our hypoexponential approximation by considering a more generic concatenation of exponentially distributed waiting times than the Erlang distribution and allowing for the rate parameters to vary between compartments. In Section~\ref{Sec:TwoMomentApprox}, we derive explicit expressions of these rates that replicates the first and second moment of the gamma distribution in \eqref{Eq:DDEIVP}. Turning to numerical results, we show that the FCRK method derived in Section~\ref{Sec:FCRK} performs to the expected accuracy in Section~\ref{Sec:NumericalVerification}. Then, leveraging the numerical simulation of the gamma distributed DDE~\eqref{Eq:DDEIVP}, we show that the hypoexponential approximation outperforms the common Erlang approximation of the underlying gamma distributed DDE in Section~\ref{Sec:ApproximationAccuracy}. As the most striking illustration, we show in Section~\ref{Sec:LinearStabilityApproximation} that the Erlang distributed DDE does not necessarily replicate the qualitative properties of the underlying gamma distributed DDE. Finally, we illustrate how to implement the hypoexponential approximation to estimate parameters of an epidemiological model in Stan in Section~\ref{Sec:StatiticalInference} before finishing with a brief discussion.

\section{Functionally continuous Runge-Kutta methods} \label{Sec:FCRK}

Most existing numerical methods for delay differential equations have been adapted from known numerical methods for ODEs \citep{Bellen2009,Eremin2017,Vermiglio1988,Enright1997}. For a given stepsize $h$ and integration mesh given by $t_n = t_0 +nh$, these continuous Runge-Kutta (RK) methods are designed to output a continuous function over the delay interval. This continuous function is then used to evaluate the solution at the abscissa $c_i$ of the RK method, which is necessary since the intermediate function evaluations in each RK step fall at precisely these points $t = t_n+c_ih - \tau.$ This illustrates another difficulty with continuous RK methods, when the delay $\tau$ is smaller than the stepsize $h$ (or vanishing, although there are additional analytical issues in this case \citep{Hartung2006,Eremin2017}),  as if $\tau < c_i h$, then overlapping will occur, i.e. the $n+1$st step will require knowledge of the solution in the current step \citep{Eremin2020,Eremin2019a}, and the method is no longer explicit. This overlapping is inevitable when solving \eqref{Eq:DDEIVP} since the convolution integral in \eqref{Eq:DDEIVP} requires knowledge of the solution $X$ on the entire semi-infinite inteval $(-\infty,t)$. To address this, continuous RK methods have been adapted to output continuous interpolants in each stage. These adapted methods are called functionally continuous Runge-Kutta methods (FCRK), and have been implemented for distributed DDEs with possibly time dependent, but finite, delay \citep{Eremin2019a,Langlois2017}. Here, we implement a 4th order FCRK method for the infinite delay problem and consider the generic initial value problem \eqref{Eq:DDEIVP}.

In what follows, we denote the function segment $X_t(\theta) = X(t+\theta) $ for $\theta <0$, and define an $s$-stage FCRK method by:
\begin{definition}[$s$-stage FCRK method]
An $s$-stage explicit FCRK method is a triple $\left( A(\alpha),b(\alpha),c\right)$ such that $A$ and $b$ are polynomial functions into $\mathbb{R}^{s \times s}$ and $\mathbb{R}^{s  }$, respectively, with $A(0) = 0$ and $b(0) = 0$, and $c \in \mathbb{R}^s$ with $c_i \geq 0.$ 
\end{definition}

It is customary to represent a $s$-stage FCRK method $\left( A(\theta),b(\theta),c\right)$ by it's Butcher tableaux
\begin{equation*}
\begin{array}
{c|cccccc}
c_i & A_{i,j}(\theta) \\
\hline 
&  b_j(\theta)  \\ 
\end{array}
\end{equation*}
where $i,j = 1,2,3,...,s $ and $A_{ij}$ and $b_j$ are the components of $A$ and $b.$ Now, for a given step size $h$, the $s$-stage FCRK method creates a continuous approximation $\eta(t) $ to the solution of the IVP~\eqref{Eq:DDEIVP} $X(t)$ through 
\begin{equation}
\eta(t) = \left \{
\begin{array}{cl}
\eta^n(h \theta) & \textrm{for} \quad t \in (t_n,t_{n+1}), \quad  \textrm{and}  \quad \theta = \frac{t-t_n}{h} \\
\psi & \textrm{for} \quad t < t_0. 
\end{array}
\right.
\label{Eq:InterpolantDefinition}
\end{equation}
The stage interpolant $\eta^n$ is a continuous approximation of the solution $X(t_n+h\theta)$ defined by
\begin{equation}
\eta^n(h\theta) = x^n + h\displaystyle \sum_{i = 1}^s b_i(\theta) K_{n,i}, \quad \theta \in (0,1) \quad \eta^0 = \psi(t_0) \quad \textrm{and} \quad x^{n} = \eta^{n-1}(h),
\label{Eq:StageStepExpression}
\end{equation}
where $ K_{i,n} = F(Y^{n,i}(c_i))$ are the stage variables and $Y^{n,i} $ is the continuous approximation of $X(t)$ in the stage given by
\begin{align*}
Y^{n,i} = x^n + h \displaystyle \sum_{j=1}^{n-1} A_{i,j}(\theta) K_{n,j}, \quad \theta \in [0,c_i].
\end{align*}

Thus, the piecewise interpolants $ \eta^n(t) $ agree with $x^n$ at the collocation points $t = \{ jh \}_{j=1}^N$ and define the piecewise $C_1$ function $\eta$. The discrete order of the numerical method is the maximal error incurred at the collocation points $jh$ 
\begin{definition}[Discrete order]
A $s$-stage method has discrete order $d$ if
\begin{align*}
\max_{t \in \{t_i\}_{i=1}^N} \| x^j - X(j h ) \| = \mathcal{O}(h^d). 
\end{align*} 
\end{definition}
Conversely, the global error of the method is the absolute error incurred throughout the simulation when considering the solution $x$ and $\eta$ as continuous functions on the interval $t \in [t_0,T]$.
\begin{definition}[Global order]
A $s$-stage method has global order $p$ if
\begin{align*}
\max_{t \in [t_0,T]} \| \eta(t) - X(t) \| = \mathcal{O}(h^p). 
\end{align*} 
\end{definition}
Now, if the local truncation error is $\mathcal{O}(h^{p+1})$, then the $s$-stage method has global order $p$ on $[t_0,T]$ \citep{Eremin2019}, and $\eta$ approximates $x_t$ to order $p$ in the sense
\begin{align*}
\max_{t \in [t_0,T]} \| \eta(t) - X(t) \| < C h^p.
\end{align*}

In what follows, we use the 4th order FCRK method  \citep{Bellen2009} with Butcher tableaux given by
\begin{equation*}
\begin{array}
{c|cccccc}
0 & 0 \\
1 & \alpha & \\
\frac{1}{2} & \alpha-1/2\alpha^2 & \frac{1}{2}\alpha^2  \\
1 & \alpha-1/2\alpha^2  &\frac{1}{2}\alpha^2  \\
\frac{1}{2} & \alpha -3/2\alpha^2 +2/3\alpha^3 & 0& 2\alpha^2 -4/3\alpha^3 & -1/2 \alpha^2 +2/3 \alpha^3 & \\ 
1 & \alpha -3/2\alpha^2 +2/3\alpha^3 & 0 & 2\alpha^2 -4/3\alpha^3 & -1/2 \alpha^2 +2/3 \alpha^3 & \\
\hline
& \alpha -3/2\alpha^2 +2/3\alpha^3 & 0 & 0 & 0 & 2\alpha^2 - 4/3\alpha^3 & -1/2 \alpha^2 + 2/3 \alpha^3  
\end{array}
\end{equation*}
although our results hold for other FCRK schemes. 

\subsection{Numerical quadrature}

In general, it is possible to adapt known FCRK methods to the infinite delay case in \eqref{Eq:DDEIVP} without significant change. However, a FCRK $s$-stage method implicitly assumes the ability to accurately calculate the right hand side of equation~\eqref{Eq:DDEIVP}. Accordingly, the main difficulty in numerically simulating \eqref{Eq:DDEIVP} is the numerical calculation of the improper convolution integral
\begin{align*}
 \int_0^{\infty} X(t-s)g_a^j(s) \d s = \int_{-\infty}^t X(s)g_a^j(t-s) \d s. 
\end{align*} 

Most numerical quadrature methods are designed for a compact domain of integration. However, artificially truncating the convolution integral in \eqref{Eq:DDEIVP} could introduce unnecessary error while simultaneously ensuring that the FCRK method is not consistent as the quadrature stepsize, $h_{\rm int}$, converges to 0. Thus, to compute the convolution integral, we map the semi-infinite domain of integration to the compact set $[0,1]$ through the the change of variables 
\begin{align*}
\omega(t,s) = \exp \left( -\frac{1}{\alpha} (t-s) ^{1/\beta} \right),
\end{align*}
where $\alpha$ and $\beta$ are two parameters determined later. The improper integral then becomes
\begin{align*}
\int_{-\infty}^t X(s)g_a^j(t-s) \d s & = \int_0^1 \frac{ \beta \alpha^{\beta j} a^j }{\Gamma(j)}  X(t- \left(-\alpha \log(\omega)\right)^{\beta} ) \exp\left[ -a \left( -\alpha \log ( \omega ) \right)^{\beta} \right] \left(-\log (\omega) \right)^{\beta j -1} \frac{\d \omega}{\omega}  \\
& = \int_0^1 u(t,\omega) \d \omega,
\end{align*}
In general, we require a $k$ times continuously differentiable integrand for a $k+1$th order Newton-Coates quadrature method to obtain it's $k+11$th order accuracy. To ensure that our change of variable does not prohibit achieving such accuracy, we show how to chose the positive constants $\alpha$ and $\beta$ to ensure that the transformed integrand is sufficiently smooth for our numerical integration techniques. 
 
\begin{lemma}\label{Lemma:AlphaBetaChoice}
Assume that $X(t)$ is $k$ times differentiable and set 
\begin{align*}
\beta = \frac{k+1}{j} +1 \quad \textrm{and} \quad \alpha = \frac{j+1}{a^{1/\beta}}.
\end{align*}
Then $u(t,\omega)$ is $k$ times differentiable in $\omega$. Further, if the $k-$th derivative of $X(t)$, $X^{(k)}(t)$, is bounded, then there exists $M$ such that
\begin{align*}
\left| \frac{\d^k}{\d\omega^k} u(t,\omega) \right| < M
\end{align*}
for $\omega \in [0,1]$.
\end{lemma}
\sloppy{The proof of the above lemma is straightforward and follows from the rapid decay of $\exp[-(-\alpha \log \left( \omega \right) ) ^{\beta} ]$  at $\omega = 0$. This decay ensures that $u(t,\omega) \to 0$ as  $\omega \to 0$. We give the full proof in Appendix~\ref{Appendix:SmoothnessConditions}. In practice, we use the $5$th order open Simpson's rule, or the $5$th order open Newton-Coates method, that requires a bounded $4$th derivative of the integrand. Therefore, when implementing Lemma~\ref{Lemma:AlphaBetaChoice}, we take $k =4$. When evaluating the numerical approximation of $I(t)$, we avoid the mesh points $t_n$ where the interpolant is continuous but not differentiable. Finally, it is known that solutions of DDEs typically have discontinuous derivatives at breaking points. However, when considering a distributed DDE such as \eqref{Eq:DDEIVP}, we can leverage the additional smoothing offered by the convolution integral and only must ensure that $t_0$ is in the integration mesh at each time point \citep{Eremin2020}.}
 
After the change of integration variable, with $\alpha$ and $\beta$ chosen as in Lemma~\ref{Lemma:AlphaBetaChoice}, solving the IVP \eqref{Eq:DDEIVP} is equivalent to solving
\begin{equation}
\left. 
\begin{aligned}
\TimeDeriv X(t) & = F \left( X(t), \int_0^{1} u(t,\omega) \d \omega \right) \\
X(s) & = \psi (s) \quad s < 0.
\end{aligned}  
\right \}
\label{Eq:TransformedGammaIVP}
\end{equation}
Then, to simulate \eqref{Eq:TransformedGammaIVP} using a FCRK method, we must numerically evaluate the convolution integral
\begin{equation}
I(t) = \int_0^{1} u(t,\omega) \d \omega.
\label{Eq:TransformedConvolutionIntegral}
\end{equation}

\subsubsection*{Quadrature rules and order conditions}

As we are developing a FCRK method to numerically integrate \eqref{Eq:TransformedGammaIVP}, we will not evaluate the transformed convolution integral \eqref{Eq:TransformedConvolutionIntegral} exactly. Rather, as mentioned, we will use a quadrature method to numerically evaluate the integral to sufficient accuracy. Specifically, we consider a FCRK method of global order $p$ so that the interpolant \eqref{Eq:InterpolantDefinition} is accurate to order $p$ on each stage. We thus have
\begin{align*}
\left| (X\ast g_a^j)(t_n) - (\eta \ast g_a^j)(t_n) \right| 
 & = \left | \int_{t_0}^{t_n} \left( X(s) - \eta(s)\right)g_a^j(t_n-s) \d s\right| \\
&\leq \left \| X(s) - \eta(s)\right \|_{L_{\infty}[t_0,T]}   \int_{t_0}^{t_n} g_a^j(t_n-s) \d s \\
& \leq Ch^p  \int_{t_0}^{t_n} g_a^j(t_n-s) \d s <  Ch^p  \int_{-\infty}^{t_n} g_a^j(t_n-s) \d s = Ch^p  
\end{align*}

Therefore, if we were to calculate the convolution integral \eqref{Eq:TransformedConvolutionIntegral} exactly, then we would evaluate the right hand side of \eqref{Eq:TransformedGammaIVP} to order $p$. In each RK stage, the evaluations of $F$ occurs within the calculation of a $K_{n,i}$, so we gain an extra order of accuracy via the multiplication by $h$ in \eqref{Eq:StageStepExpression}. Then, the local error in each step of the numerical method has order $p+1$ as required, with the extra order coming from the multiplication by $h$.

However, in practice, we cannot evaluate the convolution integral \eqref{Eq:TransformedConvolutionIntegral} exactly, but we do not need to. Rather, we use a numerical integration scheme to calculate $I(t)$. This numerical approximation must be sufficiently accurate to ensure that the overall scheme is accurate to order $p$. Indeed, as the numerical solution $\eta^n$ is only a $p-$th order approximation of the true solution $X(t)$, it is not computationally efficient to evaluate the convolution integral to extreme precision. Thus, the numerical integration should be sufficiently accurate to preserve accuracy of the method, but not so accurate as to be computationally inefficient. To illustrate this idea, assume that we evaluate the integral \eqref{Eq:TransformedConvolutionIntegral} to order $q$ using a quadrature method with stepsize $h_{\rm int}$, so
\begin{align*}
I(t) = \hat{I}(t) + \mathcal{O}(h_{\rm int}^q),
\end{align*}
 where $\hat{I}(t)$ denotes the quadrature approximation of the convolution integral. Now, consider a FCRK method of order $p$ with coefficients $K_{n,i}$ and stepsize $h$. Using Taylor's theorem, we see that
\begin{multline*}
\hat{K}_{n,1} = F(x^{n-1},\hat{I}(t_{n-1}) ) 
= F(x^{n-1}, I(t_{n-1}) + \mathcal{O}(h_{int}^q) ) 
 \\ 
 = F(x^{n-1}, I(t_{n-1})) + \partial_{x_2} F(x^{n-1},I(t_{n-1}) \mathcal{O}(h_{\rm int }^q) + \mathcal{O}(h_{\rm int}^{2q}) = K_{n,1} + \mathcal{O}(h_{\rm int}^q) \,, 
\end{multline*}
where $\partial_{x_2}F$ is the partial derivative of $F$ with respect to the second variable. So the first stage step $\hat{Y}_1$ is calculated with the same accuracy as the numerical integration. We can thus proceed inductively to calculate each $\hat{Y}_i$ and $\hat{K}_i$ with accuracy $\mathcal{O}(h_{int}^q)$. Accordingly, for the continuous approximation $\hat{\eta}^n(h\theta)$ of the solution $X(t)$, equation \eqref{Eq:StageStepExpression} gives
\begin{align*}
\hat{\eta}^n(h\theta) & = x^n + h \displaystyle \sum_{i = 1}^s b_i(\theta) \hat{K}_{n,i}  = x^n + h \displaystyle \sum_{i = 1}^s b_i(\theta) K_{n,i} + \mathcal{O}(h \times h_{\rm int}^{q}).
\end{align*}
Thus, if $h_{\rm int}^q = \xi h^p$ for some constant $\xi$, then $\mathcal{O}(h \times h_{\rm int}^{q}) = \mathcal{O}( h^{p+1})$. Therefore, the condition  $h_{\rm int}^q = \mathcal{O}(h^p)$ ensures that we do not decrease the accuracy of the scheme nor perform extra computations when numerically integrating \eqref{Eq:TransformedConvolutionIntegral}. 

Finally, we note that the integrand in \eqref{Eq:TransformedConvolutionIntegral} is not defined at $\omega = 0$. Accordingly, we use an open quadrature method so that the end points of the domain of integration, $\omega = 0$ and $\omega = 1$, are not included. In particular, we use the Simpson's open rule, also known as the 5th order open Newton-Cotes rule, which is given by
\begin{align*}
\int_a^b f(x) \d x = \frac{4h_{\rm int}}{3}\left( 2f(a+h_{\rm int}) - f(a+2h_{\rm int}) + 2f(a+3h_{\rm int}) \right) + \mathcal{O}(h_{\rm int}^5)
\end{align*}
where $h_{\rm int}  = (b-a)/4$. In practice, we use the composite rule, and must enforce $h_{\rm int}^4 = \xi h^4$, so that this numerical integration technique provides sufficient accuracy for 4th order FCRK methods. We implement this FCRK method in \citet{Matlab2017}.

\section{Ordinary differential equation approximations}

In Section~\ref{Sec:FCRK}, we developed a numerical method to solve the distributed DDE \eqref{Eq:DDEIVP}. As mentioned, numerical methods for distributed DDEs are computationally demanding, complicated and as a result, not available in most off-the-shelf scientific software packages. Therefore, we discuss a common method by which modellers avoid these difficulties via an Erlang approximation of \eqref{Eq:DDEIVP} before deriving a new phase type approximation of \eqref{Eq:DDEIVP}. 

\subsection{Erlang approximation} \label{Sec:ErlangApprox}

In many modelling applications, it is common to avoid the difficulties in simulating \eqref{Eq:DDEIVP} by enforcing that $j \in \mathbb{N}$. As previously mentioned, the case $j \in \mathbb{N}$ corresponds to $\tau^2/\sigma^2 \in \mathbb{N}.$ As $\tau^2$ being an integer multiple of $\sigma^2$ is not generic, it is common to round $j$ to the nearest integer and then set the rate parameter $b = [j]/\tau$ where $[j]$ is the nearest integer to $j$. This approximation allows modellers to replace the gamma distributed delay with an Erlang distribution and thus approximate \eqref{Eq:DDEIVP} by the Erlang distributed DDE
\begin{equation}
\left. 
\begin{aligned}
\TimeDeriv Y(t) & = F \left(Y(t), \int_0^{\infty} Y(t-s)g_b^{[j]}(s) \d s \right) \\ 
Y(s) & = \psi(s) \quad s < t_0
\end{aligned}
\right \}
\label{Eq:ErlangDDEApprox}
\end{equation}

The Erlang distributed random variable $\mathcal{E}$ with shape and rate parameters $[j]$ and $b$, respectively, has precisely the same mean as the random variable in \eqref{Eq:DDEIVP}, but not the same variance. Then, it is a simple application of the linear chain technique---where the convolution integral is written as the solution to a system of differential equations---to obtain the equivalent ODE formulation to \eqref{Eq:ErlangDDEApprox} 
\begin{equation}
\left. 
\begin{aligned}
\TimeDeriv Y(t) & = F \left(Y(t), bA_{[j]}(t) \right) \\
\TimeDeriv A_1 & = y(t)- bA_1(t) \\
\TimeDeriv A_i & = b[A_{i-1}(t)-A_i(t)] \quad \textrm{for} \quad i = 2,3,\dots,[j].
\end{aligned}
\right \}
\label{Eq:ErlangTransitCompartmentApprox}
\end{equation}

\subsection{Hypoexponenetial approximation}\label{Sec:TwoMomentApprox}

The approximation involved in the linear chain technique described previously replaces the gamma distributed convolution integral with an Erlang distributed convolution integral parameterized to match the first moment of the original gamma distribution. Here, we develop an improved approximation technique to approximate the gamma distributed DDE~\eqref{Eq:DDEIVP} by constructing a random variable $\mathcal{Y}$ with corresponding probability measure $\mathbb{P}$ that matches the first two moments of the original gamma distribution and considering the corresponding distributed DDE 
\begin{equation}
\left.
\begin{aligned}
\TimeDeriv Y(t) & = F \left(Y(t), \int_0^{\infty} Y(t-s)\d \mathbb{P}(s)\right)\\
Y(s) & = \psi (s) \quad s < t_0.  
\end{aligned}
\right \}
\label{Eq:AppoxDDEIVP}
\end{equation}
 We construct $\mathcal{Y}$ such that it represents the concatenation of exponentially distributed random variables, so it is a phase type distribution, and we show that \eqref{Eq:AppoxDDEIVP} admits a finite dimensional representation. We then derive the equivalent ODE formulation to \eqref{Eq:AppoxDDEIVP}, and show that this approximation is more accurate than the approximation in \eqref{Eq:ErlangDDEApprox}. There are infinitely many such random variables $\mathcal{Y},$ and we consider two specific cases. We discuss the benefits of each approximation in Section~\ref{Sec:SmoothedFixedComparison}.  
\subsubsection{The fixed hypoexponential approximation}

We begin by deriving the rates of the exponentially distributed random variables whose concatenation is the random variable $\mathcal{Y}_f$ where $\mathcal{Y}_f$ is the concatenation of an Erlang distribution with two exponential distributions. We parametrize the Erlang distribution so that the rates of the Erlang distribution are \textit{fixed} as the fractional part of $j$ varies. Accordingly, we refer to this approximation as the \textit{fixed hypoexponential approximation,} with corresponding random  probability measure $\mathbb{P}_f$.

\begin{theorem}\label{Thm:TwoRateApprox}
Consider the gamma distributed random variable $\mathcal{X}$ with shape parameter $j$ and mean $\tau$ and variance $\sigma^2$. Let $\mathcal{Y}_f$ be the random variable obtained by concatenating $n = \max( \ceil{j},2 )$ independent and exponentially distributed variables where $n-2$ of these exponentially distributed random variables have identical rates
\begin{align*}
\lambda_{i,f} = \frac{n}{\tau}\,,\quad i = 1,\dots, n-2\,,
%\{ \lambda_{i,f} \}_{i=1}^{n-2} =  \frac{n}{\tau}, 
\end{align*}
while the remaining two exponentially distributed variables have rates $\lambda_{n-1,f} = \nu_f$ and $ \lambda_{n,f} = \mu_f$.  Then, setting 
\begin{align*}
\nu_f = \left( \frac{\tau}{n}\left(1+\sqrt{\tfrac{n}{2j}(1-\{j\})}\right)  \right)^{-1}  
\end{align*} 
and
\begin{align*}
\mu_f = \left(  \frac{\tau}{n}\left(1-\sqrt{\tfrac{n}{2j}(1-\{j\})}\right) \right)^{-1}
\end{align*}
ensures that $\mathcal{X}$ and $\mathcal{Y}_f$ have the same first two moments.
\end{theorem}

\begin{proof}
The moment generating function (MGF) $M_{\mathcal{Y}}(t)$ of the random variable $\mathcal{Y}$ is given by 
\begin{align*}
M_{\mathcal{Y}}(t) = \prod_{i=1}^n \frac{\lambda_i}{\lambda_i-t} \quad \textrm{for} \quad t < \min_i \{\lambda_i\}.
\end{align*}
The mean $\tau_{\mathcal{Y}}$ and variance $\sigma^2_{\mathcal{Y}}$ of $\mathcal{Y}$ are therefore
\begin{align*}
M_{\mathcal{Y}}'(0) = \displaystyle \sum_{i=1}^n \frac{1}{\lambda_i} = \tau_{\mathcal{Y}} \quad \mbox{and} \quad M_{\mathcal{Y}}''(0) - \tau_{\mathcal{Y}}^2 = \displaystyle \sum_{i=1}^n \frac{1}{\lambda_i^2} = \sigma^2_{\mathcal{Y}}.
\end{align*}
 Recalling that
\begin{align*}
\lambda_i = \frac{n}{\tau}  = \lambda\,, \quad i=1,\dots,n-2
%\{ \lambda_i \}_{i=1}^{n-2} =  \frac{n}{\tau}  = \lambda,
\end{align*}
and setting  $x_1 = 1/\nu$ and $x_2 = 1/\mu$, gives
\begin{equation*} 
    x_1 + x_2  = 2\frac{\tau}{n} \quad \textrm{and} \quad  x_1^2 + x_2^2  = \frac{\tau^2}{n^2}(n(n/j-1)+2). 
\end{equation*}
From this, $x_1$ must solve 
\begin{equation*}
    x^2 - 2\frac{\tau}{n}x + \frac{\tau^2}{n^2}\bigl(1-\tfrac12 n (n/j-1)\bigr) = 0.
\end{equation*}
By symmetry, $x_2$ must be the other root of this polynomial. Hence, by denoting the fractional part of $j$ by $\{ j \} = j-\floor{j}$, we obtain 
\begin{equation*} 
    \frac{1}{\nu}  = \frac{\tau}{n}\left(1+\sqrt{\tfrac{n}{2j}(1-\{j\})}\right) \quad \textrm{and} \quad   \frac{1}{\mu}  = \frac{\tau}{n}\left(1-\sqrt{\tfrac{n}{2j}(1-\{j\})}\right), 
\end{equation*}
which ensures that the random variable $\mathcal{Y}_f$ matches the first two moments of the gamma distribution.
\end{proof}

Since the rates of the Erlang distributed random variable as $\lambda = n/\tau$ are fixed as $\{j\}$ varies, we call the random variable derived the {\em fixed hypoexponential approximation} denoted by $\mathcal{Y}_f$.  

\begin{corollary}
If the gamma distributed random variable $\mathcal{X}$ has integer shape parameter $j$, then the random variable $\mathcal{Y}_f$ defined in Theorem~\ref{Thm:TwoRateApprox} is also Erlang distributed and $\nu = \mu = \lambda_{i,f} = j/\tau$ for $i = 1,\dots,n-2$.
\end{corollary}
\begin{proof}
If $\mathcal{X}$ is Erlang distributed, then $n/j = 1$, so the unknown $x_1$ satisifies
\begin{align}
    x^2 - 2\frac{\tau}{n} x + \frac{\tau^2}{n^2} = 0.
\end{align}
Then, $x_1 = x_2 = \frac{\tau}{n},$ and the two moment approximation is exact. 
\end{proof}

\subsubsection{A smoothed hypoexponential approximation}

The parametrization of the hypoexponential distribution in Theorem~\ref{Thm:TwoRateApprox} is determined by the choice of $\{ \lambda_{i,f}\}_{i=1}^{n-2}$ and is therefore not unique. Here, we derive a slightly different parameterization of the hypoexponential approximation. This alternative approximation has benefits and a disadvantage compared to the fixed hypoexponential approximation, which we discuss below.

Again, denote the mean of the gamma-distributed random variable $\mathcal{X}$ by $\tau$ and let $j \notin \mathbb{N}$ denote the shape parameter. Now we define a second hypoexponentially-distributed random variable $\mathcal{Y}_s$ with the same mean and variance as $\mathcal{X}$. We once again use a concatenation of an Erlang distribution with two exponential distributions. Here, unlike the fixed approximation described in Theorem~\ref{Thm:TwoRateApprox}, the rate of the Erlang distribution varies continuously as the fractional part of $j$ changes. We therefore refer to this approximation as the \textit{smoothed hypoexponential approximation,}  with corresponding random  probability measure $\mathbb{P}_s$.

\begin{theorem}\label{thm:AltTwoRateApprox}
Let $\mathcal{X}$ be a $\GammaDist(j,a)$-distributed random variable where $j \notin \mathbb{N}$. Consider the hypoexponentially distributed random variable $\mathcal{Y}_s$ with rate parameters $(\lambda_{s,1}, \dots, \lambda_{s,n-2},\mu_s,\nu_s)$. Recalling that $\{j\} = j - \lfloor j \rfloor > 0$ as $j \notin \mathbb{N}$, set $\lambda_{s,1} = \dots = \lambda_{s,n-2} = j/\tau$, and define $\mu_s$ and $\nu_s$ by
 \begin{equation}\label{Eq:SmoothedRates}
    \begin{split}
    \mu_s &= \frac{2j}{\tau}\left( 1 + \{j\} + \sqrt{1 - \{j\}^2}\right)^{-1} \\
    \nu_s &= \frac{2j}{\tau}\left( 1 + \{j\} - \sqrt{1 - \{j\}^2}\right)^{-1} \\
    \end{split}
 \end{equation}
 If $j \in \mathbb{N}$, then we define $\mu_s = \nu_s =  j/\tau$. Then $ \mathcal{X} $ and  $\mathcal{Y}_s$ have the first two moments.
\end{theorem}
The proof of Theorem~\ref{thm:AltTwoRateApprox} is similar to the proof of Theorem~\ref{Thm:TwoRateApprox} and is given in Appendix~\ref{Appendix:SmoothedApproximation}.  We note that we use the term \textit{smooth} when describing the {\em smoothed hypoexpoential approximation of $X$} to refer to the continuous dependence of $\lambda_{s,i}$ on $j,$ and not in the infinitely differentiable sense. 

Once again, if $j$ is an integer, it follows from the definition that the smoothed hypoexponential approximation is exact. 

\subsection{Ordinary differential equation representation of the hypoexponential DDE}

The random variables $\mathcal{Y}_f$ and $\mathcal{Y}_s$ as defined in Theorems~\ref{Thm:TwoRateApprox} and \ref{thm:AltTwoRateApprox} correspond to the concatenation or addition of $n$ exponentially distributed random variables. As the derivation that follows is identical for the smoothed and fixed approximations, we drop the indices $f$ and $s$. The PDF of the hypoexponential distributions is obtained by convolving the PDFs of an Erlang distributed random variable with rate $\lambda$ and shape parameter $n-2$, and the two exponentially distributed random variables with respective rates $\nu$ and $\mu$, where the rates are given explicitly in Theorems~\ref{Thm:TwoRateApprox} and \ref{thm:AltTwoRateApprox}. The exponential distributions have respective PDFs $E_{\nu}$ and $E_{\mu}$. Then, the delayed term in \eqref{Eq:AppoxDDEIVP} is given by the convolution integral
\begin{align*}
\int_0^{\infty} Y(t-s) \d \mathbb{P}_i(s) = \int_0^{\infty} Y(t-s)f_{\mathbb{P}_i}(s) \d s 
\end{align*}
where $f_{\mathbb{P}_i}(s) = (g_{\lambda}^{n-2} * E_{\nu}*E_{\mu})(s)$ depends on if we are using the smoothed or fixed hypoexponential approximation. In both cases, the convolution integral 
\begin{align*}
\int_0^{\infty} Y(t-s) f_{\mathbb{P}}(s) \d s 
\end{align*}
will satisfy a system of $n$ ordinary differential equations in a similar manner to the linear chain technique \citep{Diekmann2017,Diekmann2020,Cassidy2020a}. To show that this is indeed the case, we introduce $n$ auxiliary variables $B_i(t)$ satisfying
\begin{align*}
\TimeDeriv B_1(t) & = Y(t) - \lambda B_1(t) \\
\TimeDeriv B_i(t) & = \lambda \left[B_{i-1}(t)-B_i(t) \right] \quad \textrm{for} \quad i = 2,3,..,n-2, \\
\TimeDeriv B_{n-1}(t) & = \lambda B_{n-2}(t) - \nu B_{n-1}(t) \\
\TimeDeriv B_n(t) & = \nu B_{n-1}(t) - \mu B_n(t) 
\end{align*}
with initial conditions
\begin{align*}
B_i(0) = \int_0^{\infty} \frac{\psi(-s)}{\lambda} g_{\lambda}^i(s) \d s  \quad \textrm{for} \quad i = 2,3,..,n-2
\end{align*}
and 
\begin{align*}
B_{n-1}(0) = \int_0^{\infty} \frac{\psi(-s)}{\nu} E_x(s) \d s \quad \textrm{and} \quad B_{n}(0) = \int_0^{\infty} \frac{\psi(-s)}{\mu} E_y(s) \d s.
\end{align*}

Then, using the linear chain technique on the Erlang distributed variables $B_{i}$ for $i = 1,2,3,..,,n-2$, we see
\begin{align*}
\lambda B_i(t) = (Y * g_{\lambda}^i)(t). 
\end{align*}
Then, the main result in \citet{Cassidy2020a} shows that 
\begin{align*}
\nu B_{n-1}(t) = (\lambda B_{n-2}*E_{\nu})(t) \quad \textrm{and} \quad \mu B_n(t) = (\nu B_{n-1}*E_{\mu})(t).
\end{align*}
It follows from the associativity of convolution that
\begin{align*}
\int_0^{\infty} Y(t-s)f_{\mathbb{P}}(s) \d s = (y*f)(t) = (Y*[g_{\lambda}^{n-2}*E_{\nu}*E_{\mu}])(t) = \mu B_n(t).
\end{align*}

Therefore, the distributed DDE \eqref{Eq:AppoxDDEIVP} is equivalent to the $(n+1)$ dimensional system of ODEs 
 \begin{equation}
\left.
\begin{aligned}
\TimeDeriv Y(t) & = F \left(Y(t),  \mu B_n(t) \right) \\
\TimeDeriv B_1(t) & = Y(t) - \lambda B_1(t) \\
\TimeDeriv B_i(t) & = \lambda\left[B_{i-1}(t)-B_i(t) \right] \quad \textrm{for} \quad i = 2,3,..,n-2, \\
\TimeDeriv B_{n-1}(t) & = \lambda B_{n-2}(t) - \nu B_{n-1}(t) \\
\TimeDeriv B_n(t) & = \nu B_{n-1}(t) - \mu B_n(t).
\end{aligned}
\right \}
\label{Eq:TwoMomentAppoxODEIVP}
\end{equation}
where the rates $\lambda, \mu$ and $\nu$ are taken from the fixed or smooth hypoexponential approximation.

\subsection{A comparison between fixed and smooth hypoexponential approximations}\label{Sec:SmoothedFixedComparison}

The rates $\mu$ and $\nu$ determine the expected residence time in the $n-1$st and $n$th compartments. Now, if these rates were to grow arbitrarily large, then the expected residence time would become arbitrarily small and the system of differential equations would become stiff. Furthermore, the dynamical system obtained from the gamma distributed DDE has interesting behaviour as a function of the shape parameter $j$. For $j \notin \mathbb{N},$ we expect the gamma distributed DDE to define an infinite dimensional dynamical system. However, when $j \in \mathbb{N},$ the gamma distributed DDE can be reduced to a finite dimensional system of ODEs through the linear chain technique as detailed in Section~\ref{Sec:ErlangApprox}. As $j \downarrow n-1,$ the gamma distributed DDE approaches a transit compartment model with $n-1$ compartments. However, both the fixed and smoothed approximations are equivalent to transit compartment models with $n$ compartments. Thus, it is possible that the residence time in the final compartment becomes arbitrarily small so that the extra compartment in the hypoexponential approximation is negligible at the cost of the ODE system becoming stiff. 

To formalize this argument, consider the limit of $j\to 1$ and the fixed hypoexponential distribution. Then, $\mu_f$ and $\nu_f$ must simultaneously satisfy
\begin{align*}
    \frac{1}{\mu_f} + \frac{1}{\nu_f} = \tau \quad \textrm{and} \quad \frac{1}{\mu_f^2} + \frac{1}{\nu_f^2} = \tau^2 
\end{align*}
which is only possible if  $\frac{1}{\mu_f}\times \frac{1}{\nu_f}  = 0.$ It is simple to show that, if $j>2,$ the rates $\mu_f$ and $\nu_f$ are bounded from above so that this stiffness only occurs when $j \downarrow 1$ for the fixed hypoexponential distribution. 

Now, consider the smoothed approximation and $j \downarrow n-1$ for each integer $n$. We immediately see that the rate $\mu_s$ can become arbitrary large in the limit, and the system of ODEs becomes stiff. In addition, as $j \uparrow n$, the argument of the square roots $1-\{j\}^2$ in \eqref{Eq:SmoothedRates} approaches $0$, and the derivative of $x \mapsto \sqrt{x}$ approaches $\infty$ as $x \downarrow 0$. This is problematic for optimization methods that require the gradient of the objective function. To circumvent these singularities in the smooth hypoexponential approximation, we slightly modify \eqref{Eq:SmoothedRates} by replacing $\mu_{s}$ and $\nu_s$ by $\tilde{\mu}_{s}$ and $\tilde{\nu}_s$, defined by
\begin{equation*}
\begin{split}
 \frac{1}{\tilde{\mu}_s} &= \frac{\tau}{2j} \bigl(1+\{j\} + \sqrt{1-\{j\}^2 + h^2} - \varepsilon\bigr) \\
 \frac{1}{\tilde{\nu}_s} &= \frac{\tau}{2j} \bigl(1+\{j\} - \sqrt{1-\{j\}^2 + h^2} + \varepsilon\bigr)
 \end{split}
\end{equation*}
where $0 < \varepsilon \ll 1$ is a small constant. By choosing $\varepsilon$, the practitioner can now trade-off the size of the discontinuities of the objective function at integer values of $j$, with the level of stiffness of the resulting ODEs. As we will see in Section~\ref{Sec:StatiticalInference}, for statistical inference, one often needs to optimize an objective function which depends on the solution of a DDE~\eqref{Eq:DDEIVP} at certain time points $t_0 < t_1 < t_2 < \dots < t_K$. For many optimization algorithms, it helps if the objective function depends smoothly on the model parameters, including $j$, and so using the smoothed hypoexponential in these scenarios may be advantageous.

\begin{figure}
    \centering
    \includegraphics[width=\linewidth]{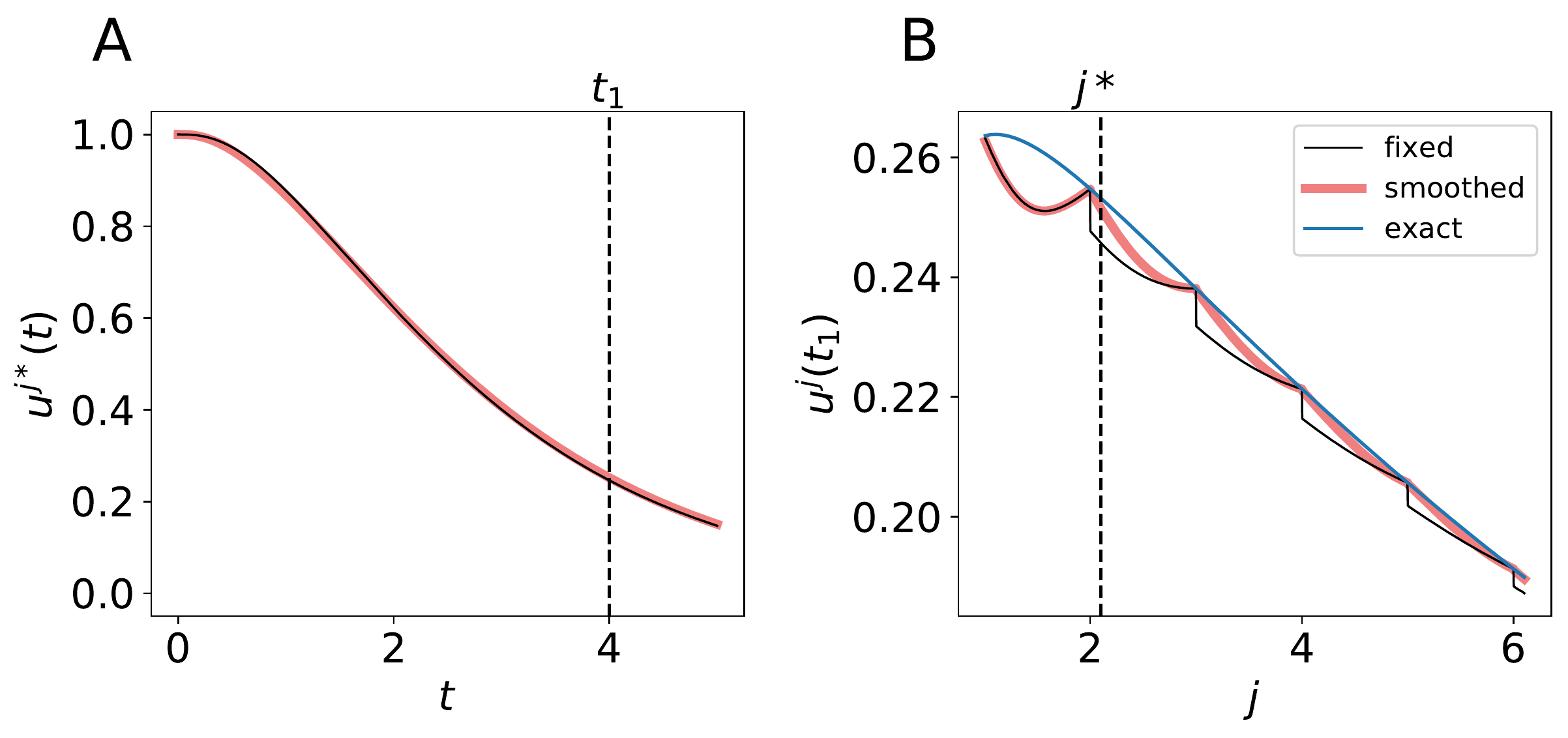}
    \caption{
        {\bf Trajectory dependence on the shape parameter $j$ is not smooth.}
        (A) The trajectory $y_f^{j^*}(t)$ and $y_s^{j^*}(t)$ for $j^* = 2.1$ calculated with the  fixed (black) and smoothed (red) hypoexponential approximations. The exact solution $u^{j^*}(t)$ is not shown as it overlaps with the two curves.
        (B) The graph of $j \mapsto u^j(t_1)$ with $t_1 = 4$, using the fixed (black) and smoothed (red) parameterizations of the hypoexponential approximation, and the exact solution (blue).
    } %% caption
    \label{fig:smooth}
\end{figure}

Further, we note that the approximations in Theorems~\ref{Thm:TwoRateApprox} and \ref{thm:AltTwoRateApprox} are approximations of the semi-infinite convolution integral in \eqref{Eq:DDEIVP}. To further compare the hypoexponential approximations, we consider the survival function of the gamma distribution with mean $1$ and shape parameter $j$ in Fig~\ref{fig:smooth} given by
\begin{align*}
    u^j(t) = \frac{j^j}{\Gamma(j)} \int_{t}^{\infty} s^{j-1}  \exp(-js) \d s,  
\end{align*}
and compute the survival functions corresponding to the fixed and smoothed hypoexponential approximations of the gamma distribution with mean $1$ and shape parameter $j$
\begin{align*}
    y_f^j(t) =    \mathbb{P}_f([t,\infty)) \quad \textrm{and} \quad y_s^j(t) = \mathbb{P}_s([t,\infty)) .
\end{align*}
We plot $y_f^j(t)$ and $y_s^j(t)$ in Fig.~\ref{fig:smooth} to illustrate the difference between the fixed and smoothed hypoexponential approximations. We do not present $u^j(t)$ as it overlaps the two approximations. Furthermore, for fixed $t = t_1,$ it is possible to view $u^j(t_1),y_f^j(t_1)$ and $y_s^j(t_1)$ as functions of $j$. In Fig.~\ref{fig:smooth} (B), we show this function for both approximations and the exact solution. For the fixed hypoexponential approximation (Theorem~\ref{Thm:TwoRateApprox}) $y_f^j(t_1)$, we generally lose continuous dependence on the parameter $j$ at integer values as the rates of the Erlang distribution $\{ \lambda_{i,f}\}_{i=1}^{n-2} = \ceil{j}/\tau$ do not vary continuously but rather jump as $j$ crosses each integer. However, using the smoothed hypoexponential parameterization, the rates $\lambda_{i,s}$ vary continuously with $j$ which appears to reduce the size of jumps at integer values of $j$. However, an analytical study of these jumps is beyond the scope of the current work.

\subsection{Approximation error estimates}

The natural phase space for distributed DDEs such as \eqref{Eq:DDEIVP} or \eqref{Eq:AppoxDDEIVP} is the space of exponentially weighted functions \citep{Cassidy2020a,Diekmann2012}
\begin{align}\label{Eq:SolutionSpace}
C_{0,\rho} = \left \{ f \in C_0 \ \bigg| \lim_{\phi \to -\infty} f(\phi)e^{\rho \phi} = 0 \right \}.  
\end{align}
In general, solutions evolving from the space of $\mathcal{X}$ measurable functions remain integrable with respect to $\mathcal{X}$ \citep{Hale1974,Cassidy2018}. 

Now, for the rate parameter of the gamma distribution given by $a = j/\tau$, solutions of the gamma distributed DDE will satisfy the growth bound in \eqref{Eq:SolutionSpace} with  $\rho < a$. Furthermore, solutions $Y(t)$ of linear gamma distributed DDE \eqref{Eq:DDEIVP} are of the form $Y(t) = Ce^{\phi t}$ with $\Re(\phi) < a$. To illustrate the increased accuracy offered by the hypoexponential approximation, we use this linear case to derive explicit bounds for the approximation error induced by replacing the gamma distribution in \eqref{Eq:DDEIVP} by an Erlang distribution as in \eqref{Eq:ErlangDDEApprox} or by the hypoexponential approximation in \eqref{Eq:AppoxDDEIVP}. In both cases, we will express the approximation error as the difference of the MGFs evaluated at $\phi$ and we will see that the hypoexponential approximation has one fewer term than the Erlang approximation. 

\subsubsection{Erlang distributed DDE}
In the Erlang approximation described in Section~\ref{Sec:ErlangApprox}, we replaced the convolution integral in the gamma distributed DDE
\begin{align*}
    \int_0^{\infty} Y(t-s) g_a^j(s) \d s 
\end{align*}
by 
\begin{align*}
   \int_0^{\infty} Y(t-s) g_b^{\ceil{j}}(s) \d s 
\end{align*}
where $b = \ceil{j}/\tau$ is the rate of the approximating Erlang distribution, and the the Erlang distributed DDE~\eqref{Eq:ErlangDDEApprox} is otherwise identical to \eqref{Eq:DDEIVP}. Thus, to compute the error induced by this approximation, we consider the difference between the convolution integrals where $Y(t) = Ce^{\phi t}$ 
% The approximation error induced by integrating with respect to the Erlang distribution in \eqref{Eq:ErlangDDEApprox} rather than the Gamma distribution is given by 
\begin{align*}
\textrm{Err}_{\mathcal{E}}(t) = \left| \int_0^{\infty} Ce^{\phi t} e^{-\phi s}\left( g_a^j(s)-g_b^{[j]}(s)\right) \d s \right|.
\end{align*}
We immediately obtain
\begin{align*}
\textrm{Err}_{\mathcal{E}}(t) =  \left| Y(t) \right| \left| M_{\mathcal{X}}(-\phi)- M_{\mathcal{E}}(-\phi) \right| 
\end{align*}
where the MGF of the gamma distributed random variable is given by 
\begin{align*}
M_{\mathcal{X}}(-\phi) = \frac{a^j}{(a+\phi)^j} = \frac{1}{(1+\phi/a)^j}, 
\end{align*}
and $M_{\mathcal{E}}(-\phi)$ is the MGF of the Erlang distribution
\begin{align*}
M_{\mathcal{E}}(-\phi) = \frac{b^{[j]}}{(b+\phi)^{[j]}} = \frac{1}{(1+\phi/b)^{[j]}}. 
\end{align*}
Then, 
\begin{align}\label{eqn:erlang-err}
\textrm{Err}_{\mathcal{E}}(t) = \left|Y(t)\right| \left| \frac{(1+\phi/b)^{[j]} - (1+\phi/a)^j}{ (1+\phi/a)^j (1+\phi/b)^{[j]} } \right|
\end{align}
Using the binomial theorem and the the fact that the Erlang distribution is parameterized so that the first moment matches that of the gamma distribution, we can write the numerator in \eqref{eqn:erlang-err} as 
\begin{align*}
    (1+\phi/b)^{[j]} - (1+\phi/a)^j 
    &= 
    \sum_{k=2}^{[j]} \left[\binom{j}{k} \left(\frac{\phi}{a}\right)^k - \binom{[j]}{k} \left(\frac{\phi}{b}\right)^k\right]  + \sum_{k=[j]+1}^{\infty} \binom{j}{k} \left(\frac{\phi}{a}\right)^k \\
    &= 
    \sum_{k=2}^{[j]} \left[ \binom{j}{k}  - \binom{[j]}{k} \left( \frac{j}{[j]} \right)^k \right] \left(\frac{\phi}{a}\right)^k + \sum_{k=[j]+1}^{\infty} \binom{j}{k} \left(\frac{\phi}{a}\right)^k
\end{align*}
Thus, the approximation error in the Erlang approximation case (see Section~\ref{Sec:ErlangApprox}) is order $(\phi/a)^2$. We see from the above analysis that if $j \in \mathbb{N}$, then \eqref{eqn:erlang-err} is identically 0 and the approximation is exact. 

\subsubsection{Hypoexponential approximations}

Turning to the two moment approximations derived in Theorems~\ref{Thm:TwoRateApprox} and \ref{thm:AltTwoRateApprox}, we see that the approximation error induced by integrating with respect to the random variable $\mathcal{Y}$ is given by 
\begin{align*}
\textrm{Err}_{\mathcal{Y}}(t) = \left| \int_0^{\infty} Ce^{\phi t} e^{-\phi s}\left[ g_a^j(s)-f_{\mathbb{P}}(s)\right] \d s \right|.
\end{align*}
 Then, we obtain
\begin{align} \notag
\textrm{Err}_{\mathcal{Y}}(t) & =   |Ce^{\phi t}| \left| \int_0^{\infty}  e^{-\phi s}\left[ g_a^j(s)-f(s)\right] \d s \right| \\
& = |Ce^{\phi t}| \left| M_{\mathcal{Y}}(-\phi) - M_{\mathcal{X}}(-\phi) \right|, \label{Eq:IntegralError}
\end{align}
where $M_{\mathcal{Y}}(-\phi)$ is the MGF of the random variable $\mathcal{Y}$ and is given by
\begin{align*}
M_{\mathcal{Y}}(-\phi) = \frac{1}{(1+\phi/\lambda)^{n-2}(1+\phi/\mu)(1+\phi/\nu)}.
\end{align*}
Then, we see 
\begin{align}\label{eqn:diff-mgf-hypo}
\left| M_{\mathcal{Y}}(-\phi) - M_{\mathcal{X}}(-\phi) \right| & =  \left| \frac{(1+\phi/a)^j - (1+\phi/\lambda)^{n-2}(1+\phi/\mu)(1+\phi/\nu)}{(1+\phi/\lambda)^{n-2}(1+\phi/\mu)(1+\phi/\nu)(1+\phi/a)^j} \right| 
\end{align}
Then, by recalling that $\phi/a < 1$ and the fact that the first two moments agree, we use the binomial theorem to write the numerator in \eqref{eqn:diff-mgf-hypo} as
\begin{multline*}
(1+\phi/a)^j - (1+\phi/\lambda)^{n-2}(1+\phi/\mu)(1+\phi/\nu) = \\
\sum_{k=3}^{n-2}  \left[ \binom{j}{k}\frac{1}{a^k} -  \binom{n-2}{k}\frac{1}{\lambda^{k}} -\left(\frac{1}{\mu}+\frac{1}{\nu}\right)  \binom{n-2}{k-1} \frac{1}{\lambda^{k-1}} - \frac{1}{\mu \nu} \binom{n-2}{k-2}\frac{1}{\lambda^{k-2}} \right]\phi^{k} 
\\
 + \left[ \binom{j}{n-1}\frac{1}{a^{n-1}} -  \binom{n-2}{n-1}\frac{1}{\lambda^{n-1}} -\left(\frac{1}{\mu}+\frac{1}{\nu}\right)  \frac{1}{\lambda^{n-2}} - \frac{1}{\mu \nu} \binom{n-2}{n-3}\frac{1}{\lambda^{n-3}} \right]\phi^{n-1} \\
 + \left[ \binom{j}{n}\frac{1}{a^{n}} - \frac{1}{\mu \nu}  \frac{1}{\lambda^{n-2}} \right]\phi^{n} + \displaystyle \sum_{k=n+1}^{\infty} \binom{j}{k}\left( \frac{\phi}{a} \right)^k . 
\end{multline*}
Now, recalling that $\lambda = n/\tau$ and $a = j/\tau$, we have 
\begin{multline*}
 \displaystyle \sum_{k=3}^{n-2}  \left[ \binom{j}{k}\frac{1}{a^k} -  \binom{n-2}{k}\frac{1}{\lambda^{k}} -\left(\frac{1}{\mu}+\frac{1}{\nu}\right)  \binom{n-2}{k-1} \frac{1}{\lambda^{k-1}} - \frac{1}{\mu \nu } \binom{n-2}{k-2}\frac{1}{\lambda^{k-2}} \right]\phi^{k}  \\
 = \displaystyle \sum_{k=3}^{n-2}  \left[ \binom{j}{k}  -  \frac{j^k}{n^k}  \binom{n-2}{k} \left( 1 + \frac{k}{n-k-1} \lambda \left( \frac{1}{\mu}+\frac{1}{\nu} \right)  + \frac{k(k-1)}{(n-k)(n-k-1)}\frac{\lambda^2}{\mu \nu}  \right) \right]\left( \frac{\phi}{a} \right)^{k}  \\
 =  \displaystyle \sum_{k=3}^{n-2}  \left[ \binom{j}{k}  -  \frac{j^k}{n^k}  \binom{n-2}{k} \left( 1 + \frac{2k}{n-k-1} + \frac{k(k-1)}{(n-k)(n-k-1)}\frac{\lambda^2}{\mu \nu}  \right) \right]\left( \frac{\phi}{a} \right)^{k}. 
\end{multline*}
Thus, as $\mu,\lambda,$ and $\nu$ are entirely determined by the mean and variance of the gamma distribution, we can write the error \eqref{Eq:IntegralError} as
\begin{align*}
\textrm{Err}_{\mathcal{Y}}(t) = \left| Y(t) \right| \left| \displaystyle \sum_{k=3}^{\infty} C_k(\tau,\sigma^2) \left( \frac{\phi}{a}\right)^k \right|
\end{align*}
where $\phi/a < 1.$ Accordingly, we see that the approximation error is order $(\phi/a )^3$, or one order better than the Erlang distributed DDE approximation. We also see that for $j \in \mathbb{N},$ as $\mu = \nu = \lambda$ as in Section~\ref{Sec:TwoMomentApprox},  $M_{\mathcal{Y}} = M_{\mathcal{X}}$ so \eqref{eqn:diff-mgf-hypo} is identically 0, and the approximation is exact.

\subsection{On three moment matching} \label{Sec:3moment}
The ODE approximations in this section aim to replicate the gamma distributed DDE by matching the first (in the case of the Erlang approximation) or first and second moments (in the hypoexponential approximations) of the underlying gamma distribution.  In Theorems~\ref{Thm:TwoRateApprox} and \ref{thm:AltTwoRateApprox}, we gave explicit expressions for the two unknown rates $\mu$ and $\nu$ to ensure that the random variables $\mathcal{Y}_f$ and $\mathcal{Y}_s$ match the first two moments of the underlying gamma distribution. For the fixed and smoothed approximations, we use four parameters,  $n = \ceil{j} = \ceil{\tau^2/\sigma^2}, \lambda_{f,s},\mu_{f,s}$ and $\nu_{f,s}$, and $n-2$ Erlang compartments having identical residence times $1/\lambda_{f,s}$. These parameters are determined by the underlying gamma distribution. We also have shown that our approximation is exact in the degenerate case where $\tau^2/\sigma^2 \in \mathbb{N}$. It is natural to ask if a similar technique could allow for a more accurate approximation by matching the first three moments. The three moment matching problem has been extensively studied \citep{Osogami2006,Bobbio2005}, and prior work indicates that it is simpler to consider the normalized first three moments given by
\begin{align*}
    m_2^N = \frac{\mathbb{E}[X^2] }{m_1^2} \quad \textrm{and} \quad m_3^N = \frac{\mathbb{E}[X^3]}{m_1 \mathbb{E}[X^2]},
\end{align*}
where 
\begin{align*}
    \mathbb{E}[X^i] = \int_0^{\infty} x^i\d  \mathbb{P}_X(x) 
\end{align*} 
denotes the expectation with respect to the random variable $\mathcal{X}$. Most existing work has centered on matching moments using an \textit{acylic phase type} distributions, of which our hypoexponential approximations are a specific type, defined by 

\begin{definition}[$n$ phase acyclic distributions]
A $n$ phase acyclic phase distribution is a phase distribution with an $n\times n$ upper triangular rate matrix.
\end{definition}

These acyclic phase distributions correspond to a Markov chain where each stage is visited at most once, i.e. the linear chain flows in one direction but some stages can be skipped. Coxian distributions are a specific type of these acyclic phase distributions that are used in modelling absorption times

\begin{definition}[Coxian distribution]
A $n$ phase Coxian distribution is an phase type distribution with an $n \times n$ bi-diagonal rate matrix $A$ that does not satisfy $a_{ii} = -a_{i,i+1}$.
\end{definition}
The Coxian distribution differs from the acyclic distrbution as each state in the corresponding Markov chain need not be visited as the absorption state can be reached from each intermediate state. We now define a $n$-phase Erlang--Coxian distribution

\begin{definition}[The $n$-phase Erlang--Coxian distribution \citep{Osogami2006}]
An $n$--phase Erlang--Coxian distribution is a convolution of an $(n-2)$ phase Erlang distribution and a two-phase Coxian distribution possibly with probability mass at zero.
\end{definition}

Let the space of distributions whose first three moments can be matched by an $n$ phase acyclic phase distribution be given by $S_n$, and let $T_n$ be the space of distributions whose first three normalized moments satisfy
\begin{enumerate}
    \item $m_2^n > \frac{n+1}{n} $ and $m_3^n \geq \frac{n+2}{(n+1)} m_2^N$
    \item $ m_2^n  = \frac{n+1}{n} $ and $m_3^n = \frac{n+2}{n}$.
\end{enumerate}
Then, \citet{Osogami2006} showed $S_n \subset T_n$.

A generic gamma distribution with shape parameter $j$ has normalized 2nd moment $(j+1)/j$. Thus, the first inequality for $m_2^n$ implies that there must be at least $n \leq j$ stages in any acyclic distribution that matches three moments of the gamma distribution. The algorithm described in \citet{Osogami2006} determines an Erlang-Coxian distribution that matches the first three moments of a generic gamma distribution with 6 free parameters and a non-zero probability of by passing the final stage in the corresponding Markov chain. In related work, \citet{Bobbio2005} constructed a 4 parameter mixed Erlang-Exponential distribution to match 3 moments of a generic gamma distribution. In the Markov chain corresponding to this mixed Erlang-Exponential distribution \citep{Bobbio2005}, there is a non-zero probability of skipping the Erlang stage and entering the absorption state immediately. 

While these algorithms admit closed form solutions in principle, they are more demanding to implement than the hypoexponential approximations derived in Theorems~\ref{Thm:TwoRateApprox} and \ref{thm:AltTwoRateApprox}. In short, their output varies depending on the ratio of normalized moments, require at least as many parameters as the hypoexponential approximations, and the non-zero probability of skipping stages does not allow for a simple skip-free Markov chain interpretation as in the hypoexponential approximation. While it may be possible to match three or more moments of a gamma distribution using a purely hypoexponential distribution with $n$ free rates, doing so requires solving increasingly large systems of polynominal equations and is therefore computationally complex and not as simple to implement as the hypoexponential approximations derived earlier. 

Along these lines, the natural extension of the hypoexponential approximations developed in Section~\ref{Sec:TwoMomentApprox} is a hypoexponential approximation where we match more than two moments of the gamma distribution. Effectively, it may appear that the $n$ phase hypoexponential approximation is easily generalized to match $m >2$ moments in the following manner: To match $m > 2$ moments of $X$, we may wish to choose $\{ x_i\}_{i=1}^m$ rates as free variables where $m <n$ and fix the final $n-m$ rates $\lambda_k = j/\tau$. In this way, we replace the system of 2 equations derived in the proof of Theorem~\ref{Thm:TwoRateApprox} by a system of $m$ polynomial equations for the $m$ moments. In Appendix~\ref{Appendix:SmoothedApproximation}, we construct a polynomial $f_m(x)$ whose roots define the unknown rates $\{ x_i\}_{i=1}^m$, and show that $f_m$ has at most two roots in $\mathbb{R}.$ This demonstrates that, when using the natural generalization of the smoothed hypoexponential distribution, it is not possible to match three or more moments of a generic gamma distribution.  

\section{Numerical results} \label{Sec:NumericalSimulation}

Here, we illustrate the analytical results of the preceeding sections and evaluate the hypoexponential approximations derived in Section~\ref{Sec:TwoMomentApprox} by comparing the direct simulation of \eqref{Eq:DDEIVP} using the FCRK method in Section~\ref{Sec:FCRK} against the numerical simulation of the approximate ODE~\eqref{Eq:AppoxDDEIVP} and the Erlang distributed DDE~\eqref{Eq:ErlangDDEApprox}. We first show that the FCRK method for \eqref{Eq:DDEIVP} is accurate to the correct order before testing the accuracy of the hypoexponential approximation derived in Section~\ref{Sec:TwoMomentApprox}.

\subsection{Numerical verification of the FCRK method}\label{Sec:NumericalVerification}

We test the 4th order FCRK numerical solver by comparing the output of the FCRK method for \eqref{Eq:DDEIVP} against differential equations with known solutions. To obtain these known solutions, we first consider \eqref{Eq:DDEIVP} in the case where the shape parameter $j$ is an integer. The gamma distribution in \eqref{Eq:DDEIVP} is thus an Erlang distribution and, using the linear chain technique, we derive an equivalent ODE formulation. This ODE formulation can either be solved analytically or simulated using established techniques for systems of ODEs as implemented in Matlab to give an expression for the reference solution $U(t)$. We also simulate the Erlang distributed DDE \eqref{Eq:DDEIVP} using our 4th order FCRK method to compute $X(t)$. Then, to compute the accuracy of our simulation, we compute the $L_{\infty}([0,T])$ error between the solution of \eqref{Eq:DDEIVP} and the solution of the equivalent ODE. In general,  a $p$th order FRCK indicates that
\begin{align*} 
E = \max_{t \in [t_0,T]} | X(t)-U(t)| \leq Ch^p,
\end{align*}
where $h$ is the stepsize of the FCRK method. The error $E$ satisfies  
\begin{align*}
\log(E) \leq \log(C) + p\log(h).
\end{align*}
Therefore, the slope of $\log(E)$ as a function of stepsize $\log(h)$ is the order $p$ of the FCRK method.

\subsubsection*{Linear test problem}

We first consider the linear test problem
\begin{equation}
\left. 
\begin{aligned}
\TimeDeriv X(t) & = \frac{4}{5}X(t)-\frac{11}{10} \int_0^{\infty} X(t-s)g_a^j(s) \d s \\
X(s) & = 1 \quad s < 0 
\end{aligned}
\right \}
\label{Eq:LinearTestProblemDDE}
\end{equation}
where we set $j \in \mathbb{N}$, and choose $a = j$ so the mean delay time $\tau = 1$. In this case, we can use the linear chain technique to reduce the Erlang distributed DDE in \eqref{Eq:LinearTestProblemDDE} to
\begin{equation}
\left. 
\begin{aligned}
\TimeDeriv U(t) & = \frac{4}{5}U(t)-\frac{11}{10} a A_j(t) \\
\TimeDeriv A_1(t) & = U(t) - a A_1(t) \\
\TimeDeriv A_n(t) & = a[A_{n-1}(t)-A_n(t)] \quad \textrm{for} \quad n = 2,3,...,j
\end{aligned}
\right \}
\label{Eq:LinearTestProblemODE}
\end{equation}
where
\begin{align*}
A_n(t) = \int_0^{\infty} \frac{Y(t-s)}{a} g_a^n(s) \d s \quad \textrm{for} \quad n = 2,3,...,j.
\end{align*}
 
Equation~\eqref{Eq:LinearTestProblemODE} is a linear system of ODEs and has an exact solution given by matrix exponentials. For $j = 1$, the analytical solution is 
\begin{align*}
X(t) = e^{-t/10} \left[ \cos \left( \frac{\sqrt{29}}{10}t \right) - \frac{2}{29} \sin \left( \frac{\sqrt{29}}{10}t \right) \right].
\end{align*}

\begin{figure} [h!]
%\noindent LinearTestFigure21Mar2021
\begin{tabular}{c}  \includegraphics[scale=0.9, trim= 50 310 00 50,clip]{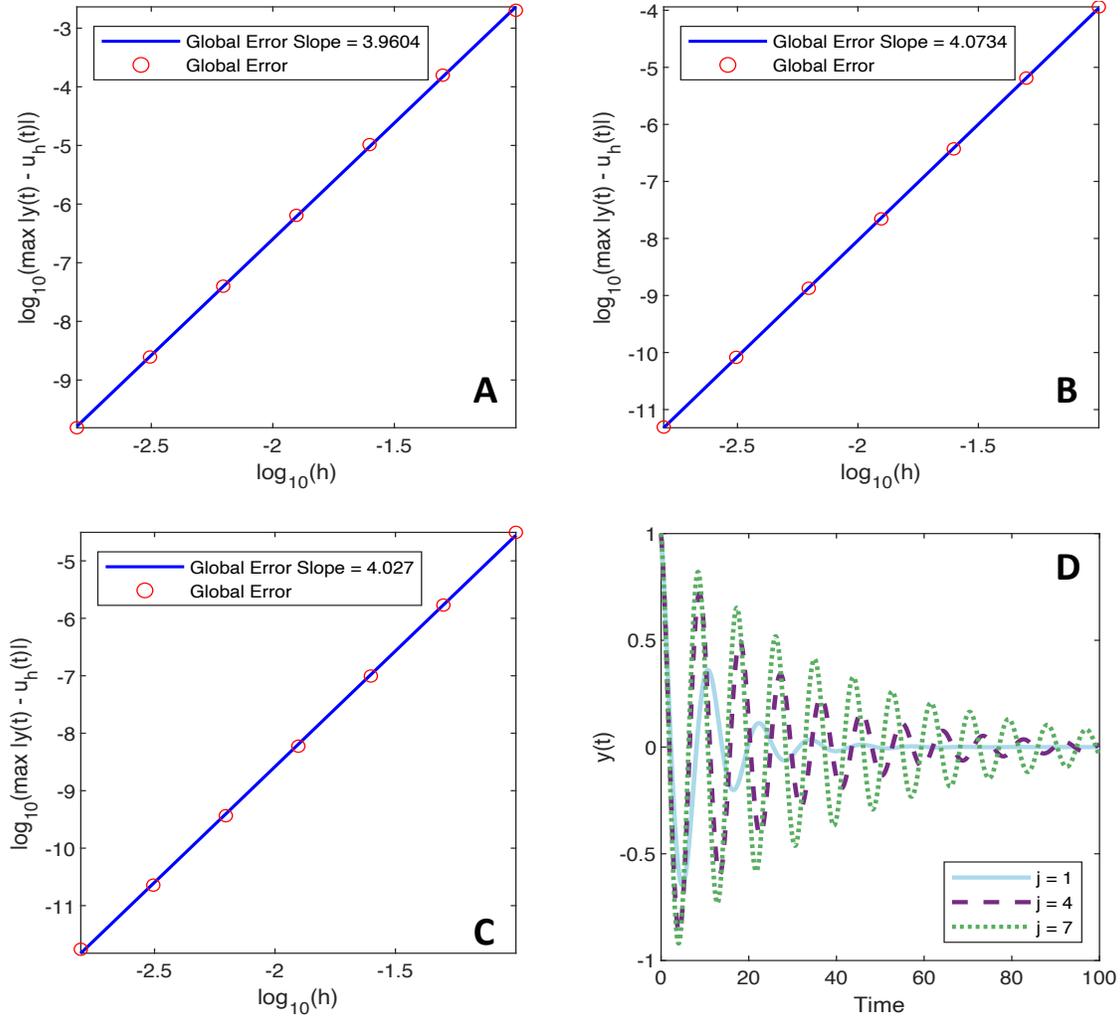} 
\end{tabular}
\caption{ Convergence plots for the linear test problem \eqref{Eq:LinearTestProblemDDE}. We plot $\log_{10}(\max_{t \in [t_0,T]} | X(t)-U(t)|)$ as a function of $\log_{10}(h)$.  The slope of $\log_{10}(\max_{t \in [t_0,T]} | X(t)-U(t)|)$ gives the convergence rate. $X(t)$ is the simulation of \eqref{Eq:LinearTestProblemDDE} using the 4th order FCRK method from Section~\ref{Sec:FCRK} with fixed step size $h$ and $U(t)$ is the solution of the equivalent ODE~\eqref{Eq:LinearTestProblemODE}. Figure (A) shows the comparison against the exact solution when $j= 1$ while figures (B) and (C) show the error between $X(t)$ and $U(t)$ for $j = 4$ and $j = 7$, respectively. Figure (D) shows the solution of the DDE for each test case. The solution $U(t)$ of the equivalent ODE~\eqref{Eq:LinearTestProblemODE} is calculated using the 4th order RK method RK45 in Matlab with relative and absolute error tolerance of $ 10^{-12}$.    }
\label{Fig:LinearDDEConvergence}
\end{figure}

Thus, we simulate \eqref{Eq:LinearTestProblemDDE} using the 4th order FCRK method described in the preceding section for $j = 1,4,7$ and compare it against the analytic solution of \eqref{Eq:LinearTestProblemODE} for $j = 1$. For $j = 4,7$, we use the 4th order variable step size RK solver in \citet{Matlab2017} with absolute and relative error tolerance of $10^{-12}$. We show the error $E(h)$ on the log-log scale and the solution of the DDE in Figure~\ref{Fig:LinearDDEConvergence}.
 
\subsubsection*{Non-linear test problem}

We next consider the non-linear test problem
\begin{equation}
\left. 
\begin{aligned}
\TimeDeriv X(t) & =  X(t)-\frac{X(t)}{K}\int_0^{\infty} X(t-s)g_a^j(s) \d s  \\
X(s) & = 1   \quad s < 0 
\end{aligned}
\right \}
\label{Eq:NonlinearTestProblemDDE}
\end{equation}
where we set $K=2$, take $j \in \mathbb{N}$, and choose $\tau= 2.25$ which gives $a = j/2.25$. Once again, we set
\begin{align*}
A_n(t) = \int_0^{\infty} \frac{Y(t-s)}{a} g_a^n(s) \d s \quad \textrm{for} \quad n = 2,3,...,j,
\end{align*}
and use the linear chain technique to reduce \eqref{Eq:NonlinearTestProblemDDE} to
\begin{equation}
\left. 
\begin{aligned}
\TimeDeriv Y(t) & = Y(t)-\frac{Y(t)(aA_n(t))}{K}  \\
\TimeDeriv A_1(t) & = Y(t) - a A_1(t) \\
\TimeDeriv A_n(t) & = a[A_{n-1}(t)-A_n(t)] \quad \textrm{for} \quad n = 2,3,...,j
\end{aligned}
\right \}
\label{Eq:NonlinearTestProblemODE}
\end{equation}

\begin{figure} [h!]
%\noindent
\begin{tabular}{c} 
 \includegraphics[scale=0.9, trim= 50 310 00 50,clip]{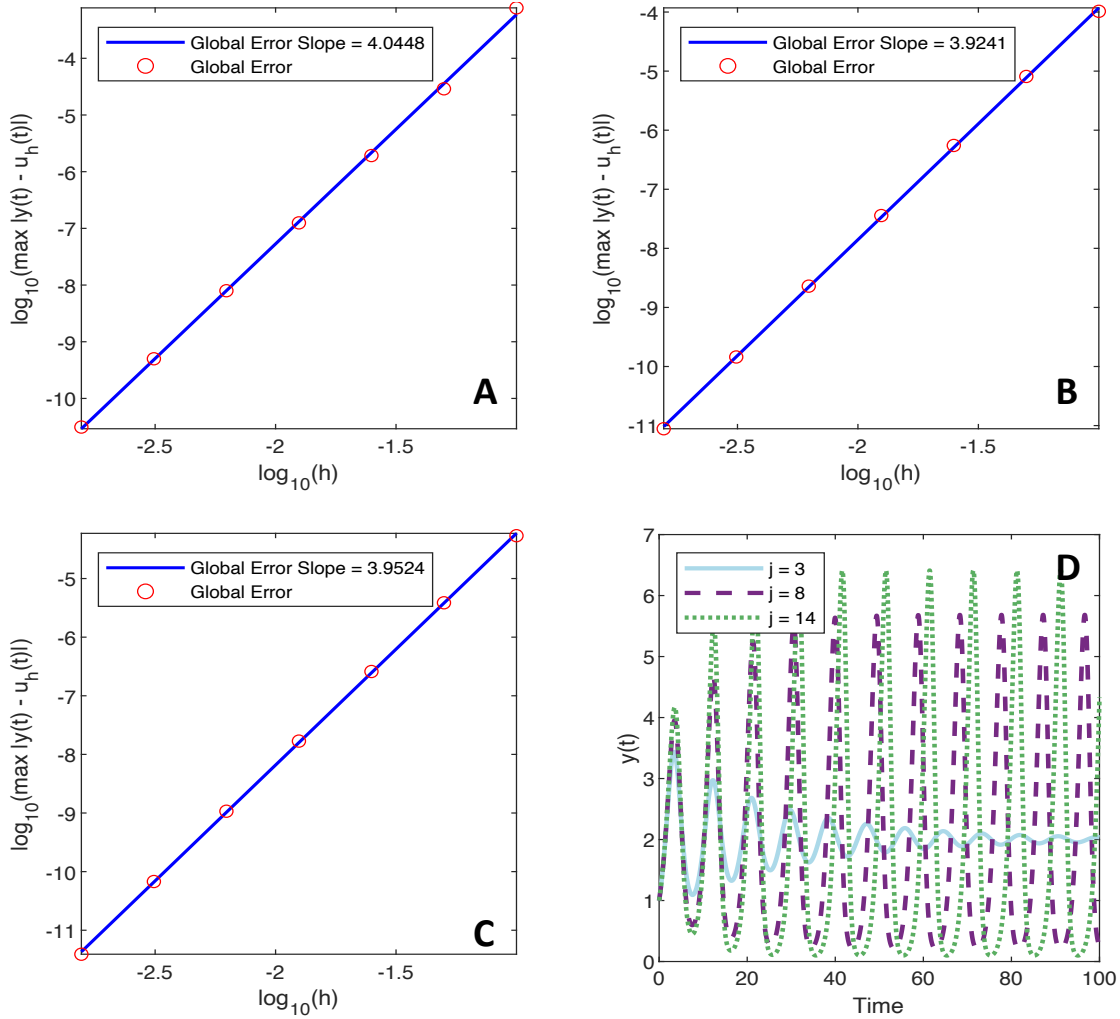} 
\end{tabular}
\caption{Convergence plots for the non-linear test problem \eqref{Eq:NonlinearTestProblemDDE}. We plot $\log_{10}(\max_{t \in [t_0,T]} | X(t)-U(t)|)$ as a function of $\log_{10}(h)$.  The slope of $\log_{10}(\max_{t \in [t_0,T]} | X(t)-U(t)|)$ gives the convergence rate. $X(t)$ is the simulation of \eqref{Eq:NonlinearTestProblemDDE} using the 4th order FCRK method from Section~\ref{Sec:FCRK} with fixed step size $h$ and $U(t)$ is the solution of the equivalent ODE~\eqref{Eq:NonlinearTestProblemODE}. Figures (A), (B) and (C) show the error between $X(t)$ and $U(t)$ for $j = 3,8,14$, respectively. Figure (D) shows the solution of the DDE for each test case. The solution $U(t)$ of the equivalent ODE~\eqref{Eq:NonlinearTestProblemODE} is calculated using the 4th order RK method RK45 in Matlab with relative and absolute error tolerance of $10^{-12}$.
}%caption
\label{Fig:NonLinearDDEConvergence}
\end{figure}

Equation~\eqref{Eq:NonlinearTestProblemODE} is a non-linear system of ODEs so we do not expect to find an analytical solution. Rather, we once again solve the system of ODEs \eqref{Eq:NonlinearTestProblemODE} using the 4th order variable step size RK solver in \citet{Matlab2017}. We solve \eqref{Eq:NonlinearTestProblemODE} with tolerance of $ 10^{-12}$, and compare this numerical solution against the numerical solution of \eqref{Eq:NonlinearTestProblemDDE} obtained using the FCRK method described in Section~\ref{Sec:FCRK}. We show the error $E(h)$ on the log-log scale for $j = 3,8,$ and $14$ and the solution of the DDE in Figure~\ref{Fig:NonLinearDDEConvergence}.

\subsubsection{Linear gamma distributed DDE}
Thus far, we have tested the FCRK method developed in section~\ref{Sec:FCRK} by simulating Erlang distributed DDEs and comparing the numerical solution against the solution of the equivalent ODE system. Here, we test our numerical method against a known solution of a gamma distributed DDE obtained by using the principle of linearised stability \citep{Diekmann2012}. In short, we consider 
\begin{align}
\TimeDeriv U(t) & = \alpha U(t) + \beta \int_0^{\infty}U(t-s)g_a^j(s) \d s.
\label{Eq:LinearGammaDistDDE}
\end{align}
We note that $U(t) = 0$ is a solution \eqref{Eq:LinearGammaDistDDE} and make the ansatz $U(t) = Ce^{\lambda t}$. Inserting $U(t)$ gives the characteristic function 
\begin{align*}
    \Delta (\lambda) = \alpha- \lambda + \beta \mathcal{L}[g_a^j](\lambda)
\end{align*}
where $\mathcal{L}[f](s)$ is the Laplace transform of the function $f$ evaluated at $s$. It follows that $\mathcal{L}[g_a^j](\lambda) = M_{\mathcal{X}}(-\lambda) $ for moment generating function of the gamma distributed random variable evaluated at $-\lambda.$ Thus, a solution of \eqref{Eq:LinearGammaDistDDE} must satisfy 
\begin{align}
     \Delta (\lambda) & = \lambda- \alpha - \beta \frac{a^j}{(a+\lambda)^j} = 0  
\end{align}
which implies 
\begin{align*}
    0 = (\alpha-\lambda)(a+\lambda)^j +\beta a^j.
\end{align*}
Now, for simplicity, we set $\alpha = -a$ so that $(a+\lambda)^{j+1} = \beta a^j$ and 
\begin{align}
\lambda = \left( \beta a^j\right)^{1/(j+1)} - a 
\label{Eq:LambdaDefinition}
\end{align}
is a solution of the characteristic function where we must impose $\beta < 2^{j+1}a $. The corresponding eigenfunction $U(t) = ce^{\lambda t}$ is the solution of the linear distributed DDE for the history function $\phi(s) = ce^{\lambda s}.$ We have thus determined an analytical solution to the linearised DDE against which we can compare the numerical solution obtained by the FCRK method.

\begin{figure} [h!]
%\noindent
\begin{tabular}{c} \includegraphics[trim= 50 590 0 0,clip,scale=0.9]{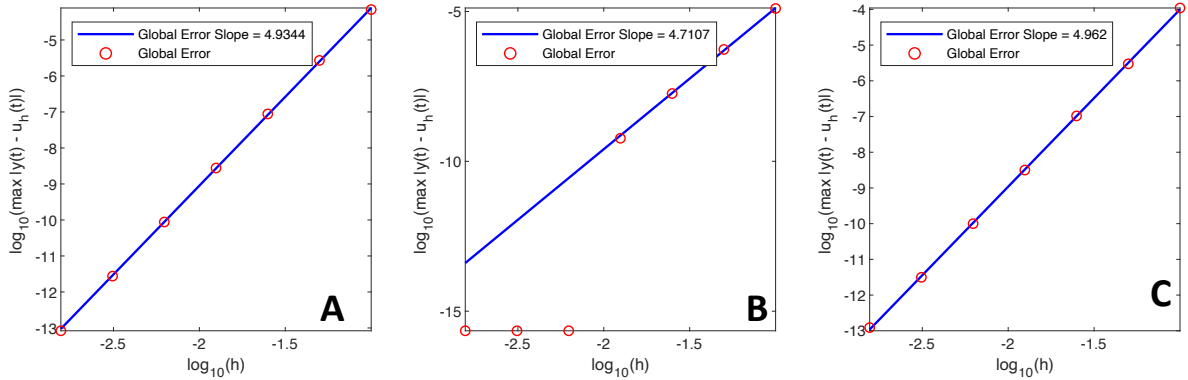} 
\end{tabular}
\caption{ Convergence plots for the linear gamma distributed DDE test problem \eqref{Eq:LinearGammaDistDDE}. We plot $\log_{10}(\max_{t \in [t_0,T]} | X(t)-U(t)|)$ as a function of $\log_{10}(h)$ where $X(t)$ is the simulation of \eqref{Eq:LinearGammaDistDDE} using the 4th order FCRK method from Section~\ref{Sec:FCRK} and $U(t)= c e^{\lambda t}$ is the analytical solution of the linear distributed DDE. In figure (B), the error reaches machine precision for $\log_{10}(h) < -2.$ Figures (A), (B) and (C) show the error between $X(t)$ and $U(t)$ for the parameter triples  $(\tau,j,\beta)$ given by   $(4.65,2.15,0.5),(3.76,3.70,0.35)$ and $(4.25,2.25,0.71)$ respectively.      }
\label{Fig:CharacteristicGammaDDEConvergence}
\end{figure}

Now, we consider parameter triples $(\tau,j,\beta)$, set $\alpha = -a = -j/\tau$ and calculate $\lambda$ by taking the principle root in \eqref{Eq:LambdaDefinition}. In Figure~\ref{Fig:CharacteristicGammaDDEConvergence}, we show the convergence of the numerical solution of \eqref{Eq:LinearGammaDistDDE} obtained using the FCRK method to the analytical solution for the parameter triples $(4.65,2.15,0.5),(3.76,3.70,0.35)$ and $(4.25,2.25,0.71)$.

We note that we observe the expected convergence rate with approximate slope $4$ in all cases until we reach numerical precision. We therefore conclude that the FCRK method derived in Section~\ref{Sec:FCRK} exhibits the expected accuracy. 

\subsection{Numerical evaluation of Erlang and hypoexponential approximations} \label{Sec:ApproximationAccuracy}

Having confirmed the accuracy of our FCRK method to solve the distributed DDE \eqref{Eq:DDEIVP}, we now evaluate the Erlang and hypoexponential approximations for the two test problems \eqref{Eq:LinearTestProblemDDE} and \eqref{Eq:NonlinearTestProblemDDE} for $j \in \mathbb{Q}$. To test the accuracy of the Erlang approximation from Section~\ref{Sec:ErlangApprox}, we use \eqref{Eq:ErlangTransitCompartmentApprox} with shape parameter $[j]$ and corresponding rate $[j]/\tau$. We also consider the fixed hypoexponential approximation as described in Section~\ref{Sec:TwoMomentApprox} with $N = \max\{ \ceil{j},2\}$ and the rates $\mu_f$ and $\nu_f$ as given in Theorem~\ref{Thm:TwoRateApprox}. In these simulations, the fixed and smoothed approximations are indistinguishable, so we only show the fixed approximation corresponding to $\mathcal{Y}_f$.

In the following simulations, we consider \eqref{Eq:LinearTestProblemDDE} and with $\tau = 1$, $j = 2.57,3.48,6.5$, and a non-constant history function given by $\phi(s) = 0.1e^{0.1s}$ for $s < 0$. We simulate the non-linear test problem \eqref{Eq:NonlinearTestProblemDDE} for $j = 2.82,4.72,$ and $6.45$ with a constant history function $\phi(s) = 0.5$. 

\begin{figure} [h!]
%\noindent
\begin{tabular}{c} \includegraphics[trim= 50 350 0 0,clip,scale=0.9]{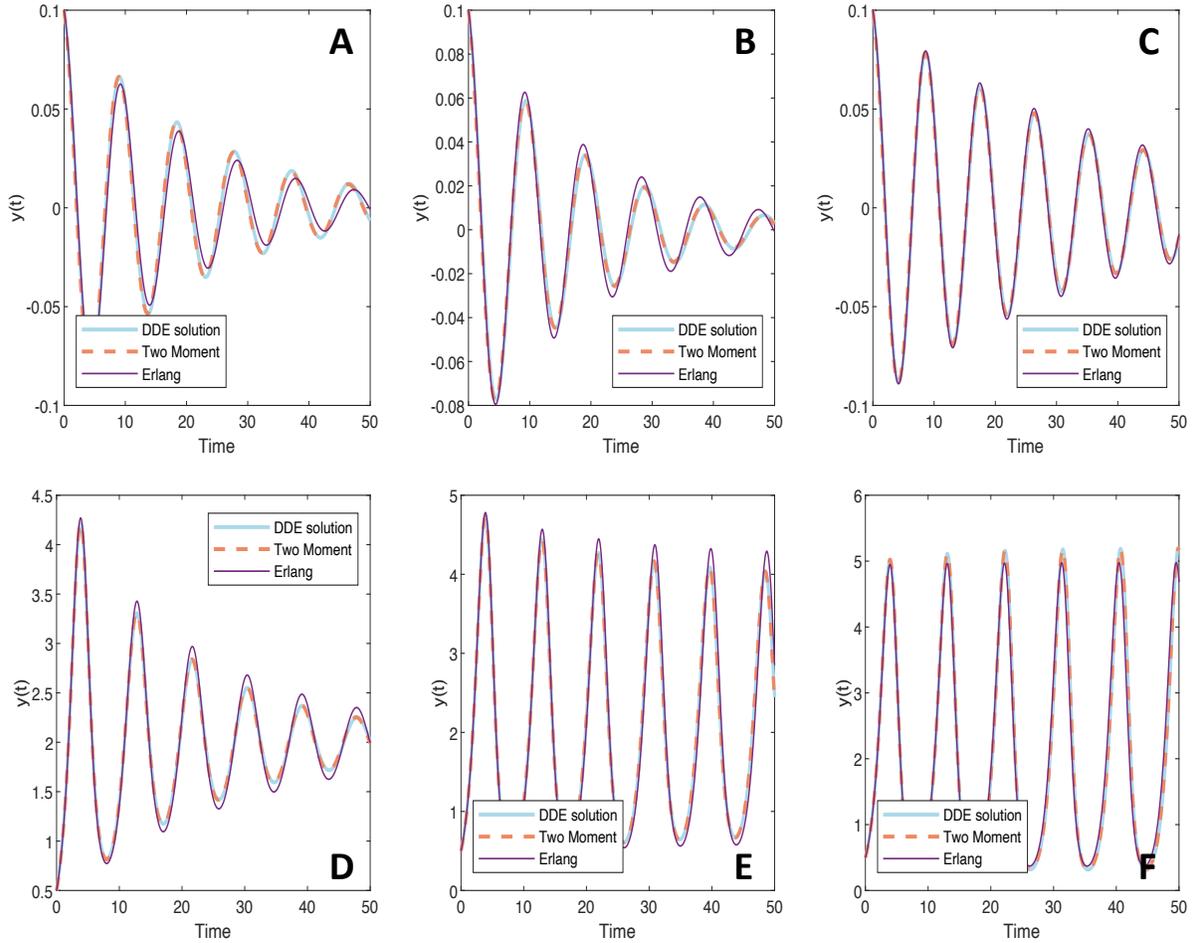} 
\end{tabular}
\caption{ Comparison of ODE approximations to the gamma distributed DDE \eqref{Eq:DDEIVP} using the Erlang approximation in equation~\eqref{Eq:ErlangTransitCompartmentApprox} or the fixed hypoexponential approximation $\mathcal{Y}_f$ in \eqref{Eq:TwoMomentAppoxODEIVP}. In all cases, the solution of the gamma distributed DDE as solved using the FCRK method is in solid blue, the solution of the fixed hypoexponential two moment approximation is in dashed orange and the solution of the Erlang approximation is in purple. Figures \textbf{A}, \textbf{B}, and \textbf{C} show the solution of the linear test problem \eqref{Eq:LinearTestProblemDDE} for $j = 2.57,3.48$ and $j = 6.5$, respectively.   Figures \textbf{D}, \textbf{E}, and \textbf{F} show the solution of the nonlinear test problem \eqref{Eq:LinearTestProblemDDE} for $j = 2.82,4.72$ and $j = 6.45$, respectively. } 
\label{Fig:ApproximationSimulationTest} 
\end{figure}

In all cases shown in Figure~\ref{Fig:ApproximationSimulationTest}, the Erlang approximation has a visbily larger error than the hypoexponential approximations. In fact, there is no perceptible difference between the fixed (and consequently, the smoothed) hypoexponential  approximation and the solution of the gamma distributed DDE. While we only present the simulation results for a limited number of test problems, this significantly improved approximation by the hypoexponential approximations was confirmed a number of other test cases. 

\subsubsection{Effects on linear stability}\label{Sec:LinearStabilityApproximation}

To study the effects of replacing the gamma distributed DDE \eqref{Eq:DDEIVP} by an Erlang or either hypoexponential approximation, we consider the linear gamma distributed DDE given in \eqref{Eq:LinearGammaDistDDE}. We note that $X(t) = 0$ is an equilibrium solution of the linear DDE, and it follows that this linear DDE represents the linearised version of
\begin{align*}
\TimeDeriv X(t) & = F \left( X(t),\int_0^{\infty}X(t-s)g_a^j(s) \d s \right).
\end{align*}
where $F(0,0)= 0$ and $\alpha = \partial_{x_1} F(x_1,x_2)|_{(0,0)}$ and $\beta = \partial_{x_2} F(x_1,x_2)|_{(0,0)}$.  The principle of linearised stability for delay equations with infinite delay was established by \citet{Diekmann2012} and, in short, indicates that the qualitative behaviour of a DDE with infinite delay near an equilibrium solution is determined by the linearised version of the DDE. 

\begin{figure} [h!]
%\noindent
\begin{tabular} {l} \includegraphics[trim= 50 550 0 0,clip,scale=0.9]{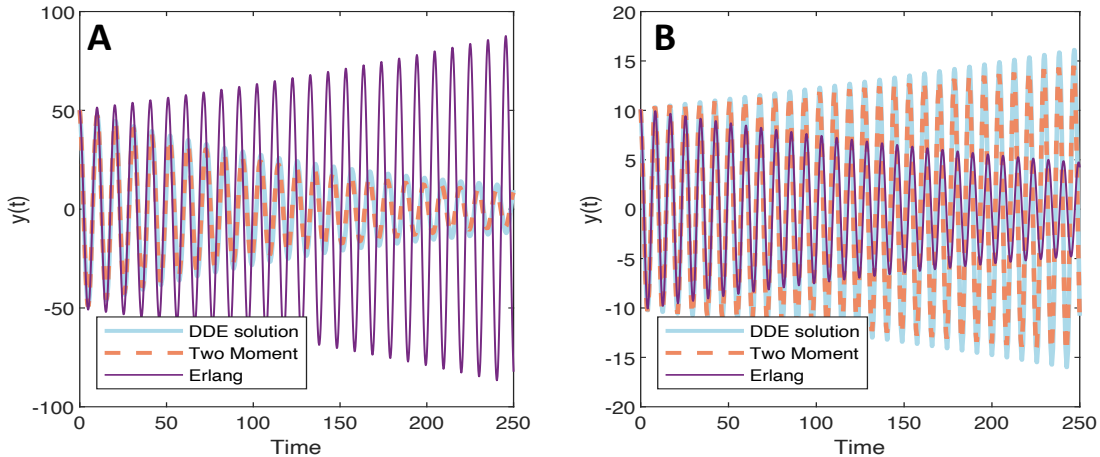} 
\end{tabular}
\caption{ Comparison of ODE approximations to the gamma distributed DDE \eqref{Eq:LinearGammaDistDDE} using the Erlang approximation in equation~\eqref{Eq:ErlangTransitCompartmentApprox} or the fixed hypoexponential two moment approximation in \eqref{Eq:TwoMomentAppoxODEIVP} showing that the Erlang approximation does not have the same stability properties as the gamma distributed DDE or the hypoexponential approximation. In all cases, the solution of the gamma distributed DDE as solved using the FCRK method is in solid blue, the solution of the hypoexponential approximation is in dashed orange and the solution of the Erlang approximation is in purple. }
\label{Fig:BifurcationApproximationTest}
\end{figure}

As a final test of the Erlang and hypoexponential approximations, we consider two specific examples with $\tau = 1$ and parameters $j,\alpha,$ and $\beta,$ and chosen near a bifurcation point. In Figure~\ref{Fig:BifurcationApproximationTest}, we show that the hypoexponential approximation has the same stability properties as the solution of the distributed DDE, but that the Erlang approximation does \textit{not} have the same stability properties. In Figure~\ref{Fig:BifurcationApproximationTest} (A), we set $j = 2.5, \alpha = 0.89$ and $\beta = -1.15$, while in Figure~\ref{Fig:BifurcationApproximationTest} (B), we set $j = 4.495,\alpha = 0.825$ and $\beta = -1.175$. We parameterize the Erlang and hypoexponential approximations as previously described in Sections~\ref{Sec:ErlangApprox} and \ref{Sec:TwoMomentApprox}.

These examples indicate that using an Erlang approximation to replace the gamma distributed DDE can introduce extreme approximation error and may not replicate the qualitative behaviour of the original gamma distributed DDE. However, these simulations also indicate that the hypoexponential approximation faithfully replicates the dynamics of the linearised gamma distributed DDE. In this sense, these results strongly advocate for the use of the hypoexponential approximation derived in Section~\ref{Sec:TwoMomentApprox} rather than the usual Erlang approximation when attempting to approximate the solution of a gamma distributed DDE. 
 
\section{Statistical inference} \label{Sec:StatiticalInference}

One benefit of the hypoexponential approximation of a gamma-distributed DDE is that it is easily implemented in existing inference software. Here we demonstrate a possible implementation, using the simple and ubiquitous example of the Kermack-McKendrick (SIR) model from epidemiology, and the probabilistic programming language Stan \citep{Carpenter2017}.

The SIR model describes compartments of susceptible ($S$), infected ($I$) and recovered ($R$) individuals, and in our version, the duration of the infectious period $T_I$ is $\GammaDist(j, j/\tau)$ distributed. Hence, in this example, we ignore any incubation time. The mean duration of the infectious period is $\mathbb{E}[T_I] = \tau$ and the variance is $\Var[T_I] = \tau^2/j$. The infection rate and the initial fraction infected are denoted $\beta$ and $\varepsilon$ respectively.
The model is then given by the following system of DDEs
\begin{equation}\label{eqn:SIR-DDE}
 \begin{split}
  \frac{\d S}{\d t} &= -\beta S I \\
  \frac{\d I}{\d t} &= \beta S I - \int_0^{\infty} \beta S(t-s)I(t-s) g^j_{j/\tau}(s)\, \d s \\
 \end{split}
\end{equation}
together with the relation $S+I+R = 1.$ We follow \citet{Champredon2018} and take initial data to kick start the epidemic at time $t = 0$ so $S(0) = 1-\varepsilon$, and $I(s) = \varepsilon \delta(s)$ where $\delta(s)$ is the Dirac delta measure at $s=0$.

Using the hypoexponential approximation of Section~\ref{Sec:TwoMomentApprox}, the above DDE model \eqref{eqn:SIR-DDE} can be replaced by the following system of ODEs
\begin{equation*}
 \begin{split}
  \frac{\d S}{\d t}   &= -\beta S I \\
  \frac{\d I_1}{\d t} &= \beta S I - \gamma_1 I_1 \\
  \frac{\d I_i}{\d t} &= \gamma_{i-1} I_{i-1} - \gamma_{i} I_i\,,\quad i = 2,\dots,n
 \end{split}
\end{equation*}
where we write $I = \sum_{i=1}^n I_i$, with $n = \lceil j \rceil$ and the rates $\gamma_i$ are given by 
\begin{equation}\label{eqn:defn-rates}
\frac{1}{\gamma_i} = \left\{\begin{array}{ll} 
\tau / n  & \mbox{if } i \leq n-2 \\
\frac{\tau}{n}\left(1 + \sqrt{\tfrac{n}{2j}(1-\{j\}}\right) & \mbox{if } i = n-1 \\
\frac{\tau}{n}\left(1 - \sqrt{\tfrac{n}{2j}(1-\{j\}}\right) & \mbox{if } i = n \\
\end{array}\right.
\end{equation}
where $\{j\} = j - \lfloor j\rfloor$ denotes the fractional part of $j$.
As initial condition, we take $S(0) = 1-\varepsilon$, $I_1(0) = \varepsilon$, and
$I_i(0) = 0$ for $i=2,\dots, n$.

\begin{figure}[ht!]
 \includegraphics[width=\linewidth]{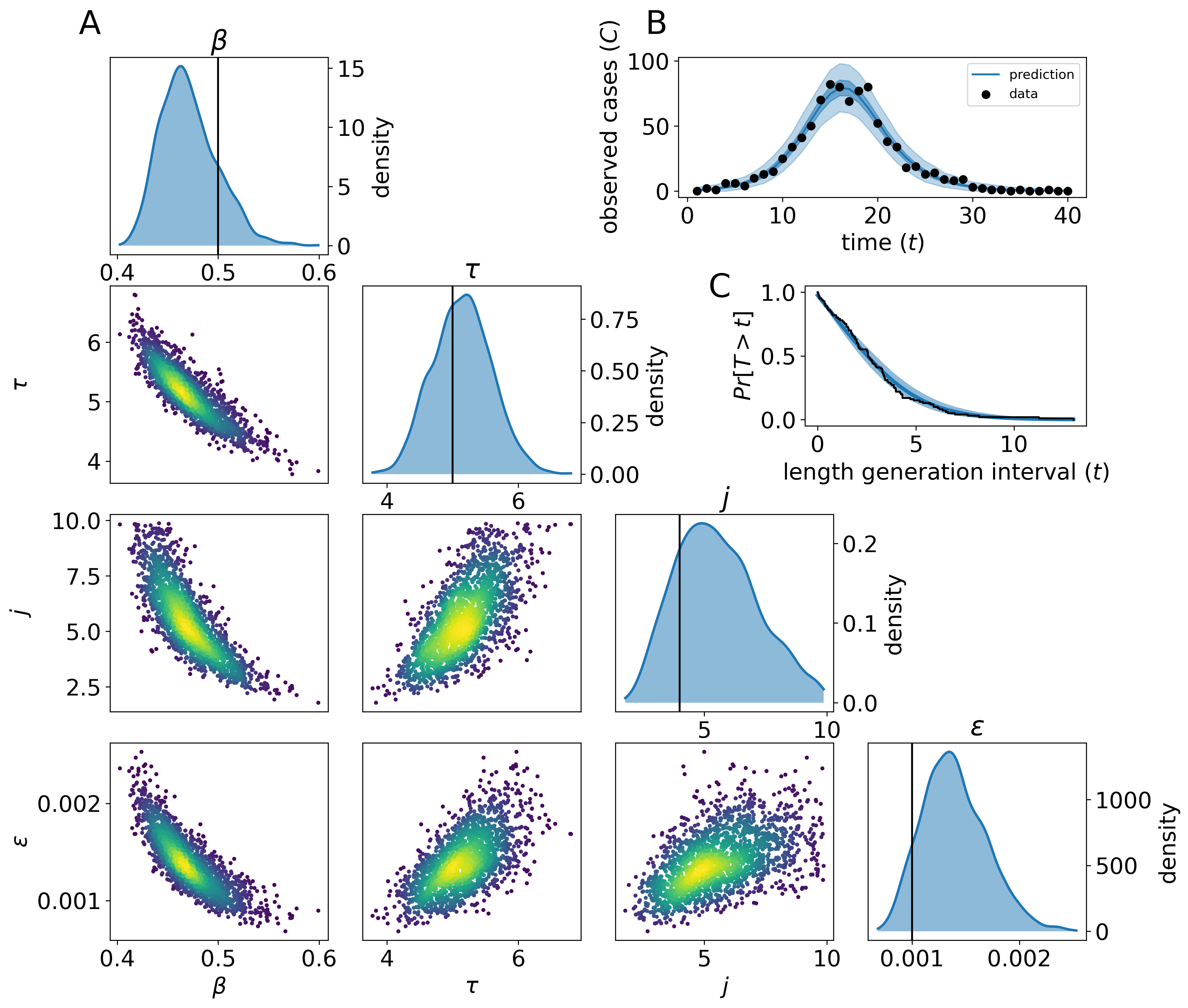}
 \caption{\label{fig:fit-posterior}
 	{\bf SIR model fit and parameter estimates.}
 	(A)	The marginal and joint posterior density of the SIR model
 	parameters.
    Each dot in the joint density scatter plots represents a
    Monte-Carlo sample from the posterior distribution. The color 
    of the dots indicates the density. The black vertical lines
    represent the ground-truth parameter values.
    (B) Simulated data $C_k$ and the model prediction $\Delta S_k M$.
    The dark-blue band represents the $95\%$ credible interval,
    and the light-blue band the $95\%$ prediction interval.
    (C) Simulated serial intervals represented as a empirical 
    survival function (black), and the fitted survival function corresponding to
    $h_{j/\tau}^j(t)$ (blue).
    The ground-truth parameters used for the simulation are
    $j=4$, $\beta = 0.5$, $\epsilon = 10^{-3}$, $\tau = 5$, $M = 10^3$, 
    and $L = 10^2$.
 }%%caption
\end{figure}

Using the above model \eqref{eqn:SIR-DDE}, we use the cumulative incidence 
$\Delta S_k \equiv S(t_{k-1}) - S(t_k)$ to simulate cases
\begin{equation*}
    C_k \sim \Poisson(\Delta S_k M)
\end{equation*} 
where $M$ is a large constant and $t_1 < t_2 < \dots < t_K$ are (positive) observation times. The simulated incidence data is shown in Fig~\ref{fig:fit-posterior}B.

In addition to time series of the number of reported cases,
often other data is collected to inform an epidemic model.
For example, symptom onset data from transmission couples 
might be available, which gives information about the 
length of the generation interval $T_G$.
Assuming that the duration of the infection is gamma distributed,
the hypoexponential approximation method allows one to estimate
the shape parameter of this distribution using both 
time series data and transmission couple data simultaneously,
i.e.\ using ``evidence synthesis''.

Suppose that we have $L$ observed serial intervals. Assuming for simplicity that symptom onset is immediate, the distribution of the 
serial interval and generation interval are identical and the 
density function for $T_G$ is given by \citep{Svensson2007}
\begin{equation*} 
 h_{j/\tau}^j(t) = \frac{1}{\tau} \int_t^{\infty} g_{j/\tau}^j(s) \d s
\end{equation*}
The log-likelihood of the data is now the sum of the log-likelihood of the case data $C = (C_1,\dots, C_K)$, and the log-likelihood of the transmission-couple data $T = (T_1,\dots, T_L)$ 
\begin{equation*}
 \mathcal{L}\left[C, T | \theta \right] = \sum_{k=1}^{K} \log p_{M \cdot \Delta S_k }(C_k) + \sum_{\ell=1}^{L} \log h_{j/\tau}^j(T_{\ell})
\end{equation*}
where $p_{\mu}(x) = e^{-\mu} \mu^x/x!$ is the density of the Poisson distribution.
To demonstrate this approach, we simulated, in addition to observed cases $C$, a small number of generation times $T_{\ell} \sim h_{j/\tau}^j$.
The simulated serial interval data is shown in Fig~\ref{fig:fit-posterior}C.

We then used Stan to fit the model to the simulated data in a Bayesian framework. The fitted model predictions are shown together with the simulated data in Fig~\ref{fig:fit-posterior}B and C.
The marginal and joint posterior densities of the model parameters 
are shown in Fig~\ref{fig:fit-posterior}A, together with the ground-truth 
values used to simulate the data.

The Stan model code and a python script to simulate data and fit the 
model is available on \url{https://github.com/lanl/gamma-dde}.

\section{Discussion}

Gamma distributed DDEs, such as \eqref{Eq:DDEIVP}, occur throughout mathematical biology. However, modellers often make simplifying assumptions due to the lack of appropriate numerical methods for infinite delay models. In this work, we developed a FCRK method to numerically simulate gamma distributed DDEs, established order conditions to ensure accuracy of the method, and illustrated our results with a series of test problems. Despite the development of a FCRK method to simulate \eqref{Eq:DDEIVP} in this work, many software packages rely on ODE solvers to perform parameter fitting and statistical inference. Accordingly, we derived a finite dimensional approximation of the gamma distributed DDE using a hypoexponential approximation and used numerical simulation to show that this hypoexponential approximation outperforms the common Erlang approximation. In particular, we demonstrated that using the Erlang approximation can lead to qualitatively different behaviour than the hypoexponential approximation and the true solution of the gamma distributed DDE. Finally, we implemented our finite dimensional approximation in Stan \citep{Carpenter2017} to fit synthetic data from a hypothetical epidemic.

The primary impediment towards the adaptation of FCRK method to distributed DDEs with infinite delay is the accurate and consistent evaluation of the convolution integral \eqref{Eq:TransformedConvolutionIntegral}. Here, we developed a change of variable that transforms the semi-infinite domain of integration to $[0,1].$ This change of variables is parametrized by certain parameters $\alpha$ and $\beta$, and we give conditions on $\alpha$ and $\beta$ on to ensure that the transformed integrand is sufficiently smooth to implement standard quadrature rules in Lemma~\ref{Lemma:AlphaBetaChoice}. These conditions, and thus the change of variable and FCRK method, apply to other delay distributions that decay exponentially. Accordingly, the FCRK method framework developed in this article should extend to other distributed DDEs with infinite delays with minimal changes. 

Until such FCRK methods are implemented in common software packages such as Stan, it is useful to have accurate finite dimensional approximations of the distributed DDE \eqref{Eq:DDEIVP}. Many finite dimensional approximations have been developed recently. The hypoexponential approximation described in this work offers a number of advantages over existing methods. While our analysis of the approximation error does not allow for an explicit expression of the error introduced by replacing the gamma distribution by either the Erlang or hypoexponential distribution, it offers a heuristic explanation for why the hypoexponential distribution is more accurate than the common Erlang approximation. As we show via numerical simulation, the hypoexponential approximation does not result in solutions that are qualitatively different from the underlying gamma distributed DDE for our test problems, unlike the common Erlang approximation. Moreover, unlike existing algorithms to parametrise phase type distributions that do not give explicit values for the parameters of the phase type distribution, we explicitly derived the parameters of the hypoexponential distribution as a function of the mean and variance of the underlying gamma distribution. This explicit expression for the rates allows for simple implementation in software packages such as Stan and we showed how to implement a simple SIR model with a gamma distributed duration of infection. 

In the SIR model, the shape parameter $j$ of the infectious period $T_I$ is hard to identify due to the correlation with other parameters, as shown in Figure~\ref{fig:fit-posterior}. Also, the trajectories of the model as a function of $j$ are very similar when $j$ is large, which makes the likelihood of the incidence data $C$ not well behaved. This problem has also been described for {\it in-vitro} SHIV data \citep{Beauchemin2017}.  To resolve such identifiability issues, it can be important to use other data to inform the parameter $j$. In our example, we used synthetic serial intervals that could be observed during real-life epidemics using transmission pairs. The distribution of these serial intervals depends on the real-valued parameter $j$. Therefore, if one wants to fit the model to both incidence data and serial intervals in an evidence synthesis framework, it is essential that the likelihood of the incidence data also depends on a real-valued shape parameter~$j$. The hypoexponential distribution therefore facilitates such evidence synthesis.

Treating $j$ as a first-class real-valued parameter in a statistical model can be important for accurately and efficiently estimating important quantities as the basic reproduction number $R_0$. Furthermore, for certain childhood diseases it can be shown that the dependence on $j$ is of a more qualitative nature \citep{Krylova2013} due to bifurcations. As we have shown, using an Erlang instead of the hypoexponential approximation in such cases, can result in large deviations from the true gamma distributed model.

Our work represents a step towards the relaxing the assumption of Erlang distributed delays in, amongst many other applications, infectious disease epidemiology. The FCRK method developed in this work allows modellers to directly simulate a gamma distributed DDE if precise numerical results are necessary, while the hypoexponential approximations offer a more accurate ODE representation of the underlying DDE than the common Erlang approximation without any increase in complexity. Accordingly, we have presented two distinct pathways to allow for the implementation of gamma distributed DDEs and facilitate their use in mathematical biology or other fields of science.

\section*{Acknowledgements}
Portions of this work were done under the auspices of the U.S. Department of Energy under contract 89233218CNA000001 and supported by National Institutes of Health (www.nih.gov) grants R01-OD011095 (CHvD, TC) and R01-AI116868 (TC). PG was supported by a National Science and Engineering Research Council (NSERC) Undergraduate Student Research Award. ARH was funded by NSERC Discovery Grant RGPIN-2018-05062. 

% \clearpage
%
%\bibliography{library,addl-refs}

\appendix
\section{Smoothness conditions for the FCRK}\label{Appendix:SmoothnessConditions}
Here, we give the proof of Lemma~\ref{Lemma:AlphaBetaChoice} where we gave sufficient conditions on the change of variable
    \begin{align*}
\omega(t,s) = \exp \left( -\frac{1}{\alpha} (t-s) ^{1/\beta} \right)
\end{align*}
to ensure that the transformed convolution integral is sufficiently smooth to not introduce unnecessary error in the FCRK. We recall that $X(t-s)$ is the solution of the gamma distributed DDE and $g_a^j(s)$ is the PDF of the gamma distribution with shape parameter $j$ and rate parameter $a$, and we are calculating
\begin{align*}
      \int_0^1 \frac{ \beta \alpha^{\beta j} a^j }{\Gamma(j)}  X(t- \left(-\alpha \log(\omega)\right)^{\beta} ) \exp\left[ -a \left( -\alpha \log ( \omega ) \right)^{\beta} \right] \left(-\log (\omega) \right)^{\beta j -1} \frac{1}{\omega} \d \omega  = \int_0^1 h(t,\omega) \d \omega
\end{align*}

We first show that the derivatives of $h(t,\omega)$ with respect to $\omega$ can be computed inductively.

\begin{lemma}
 Assume that $X(t)$ is $k$ times differentiable. Then,
 \begin{align*}
     \frac{ \d^n }{\d \omega^n} h(t,\omega) = \frac{\beta a^j\alpha^{\beta j} }{\Gamma(j) } \displaystyle \sum_{i=0}^{4^n} C_{n_i} X^{d_{n_i}}\left (t-\left[ -\alpha\log(\omega) \right]^{\beta} \right) \frac{\exp\left[-a (-\alpha\log(\omega))^{\beta}\right] (-\log(\omega) )^{\beta j-(n+1) +b_{n_i} } }{ \omega^{n+1} }
 \end{align*}
 where $C_{n_1}$ is a constant depending only on $n_i$, $\beta >1, b_{n_i} \geq 0, \alpha > 0$ and $d_{n_i}$ is an integer between $0$ and $n$, inclusive. 
\end{lemma}

\begin{proof}
    The proof is by induction on the order of the derivative $n$. The $n=0$ case follows immediately from the definition of $h(t,\omega)$, while the $n+1$st case comes from term by term differentiation with
    \begin{equation*}
    \begin{array}{llll}
         C_{(n+1)_{4i-3}} = C_{n_i}\alpha^{\beta}\beta;  &  C_{(n+1)_{4i-2}} = C_{n_i} a \alpha^{\beta} \beta; &   C_{(n+1)_{4i-1}} = C_{n_i}(\beta j -(k+1) + b_{k_i};  & C_{(n+1)_{4i} } = -C_{n_i} \\ 
     b_{(n+1)_{4i-3}} = b_{n_{i}}+\beta;  & b_{(n+1)_{4i-2}} = b_{n_{i}}+\beta ; & b_{(n+1)_{4i-1}} = b_{n_{i}} ;&   b_{(n+1)_{4i}} = b_{n_{i}} +1 ;\\
      d_{(n+1)_{4i-3}} = d_{(n)_{i}} +1; &   d_{(n+1)_{4i-2}} = d_{n_{i}}; &  d_{(n+1)_{4i-1}} = d_{n_{i}}; & d_{(n+1)_{4i}} = d_{n_{i}}. \\
     \end{array}
    \end{equation*}
\end{proof}

We now must show that $h(t,\omega)$ is a bounded function of $\omega$ on $[0,1]$. Take $\epsilon >0$ and consider the compact interval $\Omega$ away from $0$, $ \Omega = [\epsilon,1].$ We note that $h(t,\omega)$ is a product of continuous functions, and thus continuous, so the image of a compact set is compact and thus bounded. Then, we consider the interval $[0,\epsilon),$ and define
\begin{align*}
    \xi(x) = \exp \left( -a [-\alpha \log(x)]^{\beta} \right) \left(-\alpha \log(x)\right)^{\gamma} \frac{1}{x^{\delta} },
\end{align*}
for $\delta, \gamma >0$. We note that, for $\delta = n+1$ and $\gamma = \beta j - (n+1) +b_{n_i}$, $\xi(x)$ appears in the derivative of $h(t,\omega)$. Furthermore, since we are assuming that the solution of the DDE, $X(t)$, is continuously differentiable, $\xi(x)$ is the only possible term that can lead to blow up of $h(t,\omega)$ or it's derivatives. Now, $\xi(x) > 0 $ for $x\in(0,1)$, and we compute
\begin{align*}
    \lim \limits_{x\to 0}  -\alpha x\log(x) = 0
\end{align*}
by l'H\^{o}pital's rule. Thus, for $x$ sufficiently close to $0$ and $\gamma >0$, we have $\left( -\alpha \log(x) \right) ^{\gamma}  < 1/x^{\gamma} $ and we can therefore bound $\xi(x)$ from above as
\begin{align*}
    \xi(x) \leq \exp \left( -a [-\alpha \log(x)]^{\beta} \right)  \frac{1}{x^{\delta+\gamma} } = \exp \left[ -\left(\delta+\gamma \right) \log(x) -a [-\alpha \log(x) ]^{\beta} \right]  
\end{align*}
 Then, as $x \to 0,$ we have $\nu = -\alpha \log(x) \to \infty,$ and we see
 \begin{align*}
     0 \leq  \lim \limits_{x\to 0} \xi(x) \leq \lim \limits _{\nu \to \infty} \exp\left[ \frac{\delta+\gamma}{\alpha} \nu - a \nu^{\beta} \right] = 0  
 \end{align*}
as $\beta > 1.$ It follows that $\xi(x)$ is bounded on the entire interval $[0,1].$ The condition $\gamma >0$ is crucial in the above calculation, as it implies that $\beta j - (n+1) +b_{n_i} > 0,$ and leads to the result. 

\begin{lemma}\label{Lemma:AlphaBetaChoiceAppendix}
Assume that $X(t)$ is $k$ times differentiable and set 
\begin{align*}
\beta = \frac{k+1}{j} +1 \quad \textrm{and} \quad \alpha = \frac{j+1}{a^{1/\beta}}.
\end{align*}
Then $h(t,\omega)$ is $k$ times differentiable in $\omega$. Further, if the $k-$th derivative of $X(t)$, $X^{(k)}(t)$, is bounded, then there exists $M$ such that
\begin{align*}
\left| \frac{\d^k}{\d\omega^k} h(t,\omega) \right| < M
\end{align*}
for $\omega \in [0,1]$.
\end{lemma}

\begin{proof}
 We recall that $b_{n_i} \geq 0,$ so taking $\beta = \frac{k+1}{j} + 1 $ ensures that $\beta >1 $ and $\beta j - (n+1) +b_{n_i} > 0.$ The bound of $h^{(k)}$ follows from the boundedness of $\xi(x)$ demonstrated previously. 
\end{proof}

\section{The smoothed hypoexponential approximation} \label{Appendix:SmoothedApproximation}

Here we give the proof of Theorem~\ref{thm:AltTwoRateApprox}. As in the main text, let $\tau$ denote the mean of a Gamma distributed random variable $X$, and let $j$ denote the shape parameter, such that $X$ has variance $\sigma^2 = \tau^2 / j$. We write $a = j/\tau$ for the rate parameter, and $n = \lceil j \rceil$ for the smallest integer greater than $j$. We recall the definition of the smoothed hypoexponentially distributed random variable $Y_s$ with the same mean and variance as $X$

\begin{theorem}\label{thm:AltTwoRateApproxAppendix}
Let $\mathcal{X}$ be a $\GammaDist(j,a)$-distributed random variable where $j \notin \mathbb{N}.$ Consider the hypoexponentially distributed random variable $\mathcal{Y}_s$ with rate parameters $(\lambda_{s,1}, \dots, \lambda_{s,n-2},\mu_s,\nu_s)$. Recall that $\{j\} = j - \lfloor j \rfloor > 0$ as $j \notin \mathbb{N}$, set $\lambda_{s,1} = \dots = \lambda_{s,n-2} = j/\tau$, and define $\mu_s$ and $\nu_s$ by
 \begin{equation}\label{Eq:SmoothedRatesAppendix}
    \begin{split}
    \frac{1}{\mu_s} &= \frac{\tau}{2j}\left( 1 + \{j\} + \sqrt{1 - \{j\}^2}\right) \\
    \frac{1}{\nu_s} &= \frac{\tau}{2j}\left( 1 + \{j\} - \sqrt{1 - \{j\}^2}\right) \\
    \end{split}
 \end{equation}
 Then  $\mathcal{X}$ and $\mathcal{Y}_s$ have the same first two moments.
\end{theorem}

\begin{proof}[Proof of Theorem~\ref{thm:AltTwoRateApproxAppendix}]
    The mean of $X$ and $Y_s$ are given by
    $\mathbb{E}[X] = \tau$ and 
    \begin{equation*}
        \mathbb{E}[Y_s] = \kappa_1^Y = \sum_{k=1}^n a_k^{-1} = (n-2) \frac{\tau}{j} + 2 \cdot \frac{\tau}{2j}(1 + \{j\}) = \frac{\tau}{j}(n-1 + \{j\})
    \end{equation*}
    because the square-roots in Eq~\eqref{Eq:SmoothedRatesAppendix} cancel.
    Now using the fact that $n = j - \{j\} + 1$, we indeed find that $\mathbb{E}[Y] = \tau$.
    The variance of $X$ is equal to $\Var[X] = \tau^2 / j$ and the variance of $Y_s$ is given by
    \begin{equation*}
        \Var[Y_s] = \kappa_2^Y = (n-2) \frac{\tau^2}{j^2} + \frac{1}{a_{n-1}^2} + \frac{1}{a_{n}^2}
    \end{equation*}
    For any two numbers $u$ and $v$, we have $(u+v)^2 + (u-v)^2 = 2(u^2 + v^2)$. Hence, we find that
    \begin{equation*}
        \frac{1}{\mu_s^2} + \frac{1}{\nu_s^2} = \frac14 \frac{\tau^2}{j^2} 2\left((1+ \{j\})^2 + 1-\{j\}^2\right) = \frac{\tau^2}{j^2}(1 + \{j\})
    \end{equation*}
    Therefore $\Var[Y] = \tau^2/j^2 \cdot (n-2 + 1 + \{j\}) = \tau^2 / j$, which proves the theorem.
\end{proof}

\begin{remark}
Notice that the values $\mu^{-1}$ and $\nu^{-1}$ are the roots of the quadratic polynomial 
\begin{equation*}
    x^2 - x \frac{\tau}{j}(1 + \{j\})  + \frac{\tau^2}{j^2} \frac12 (1 + \{j\}) \{j\}.
\end{equation*}
\end{remark}
 
\subsection{Matching more than two moments}

A natural extension suggested by the hypoexponential approximations developed above, is an approximation where we match more than two moments with the Gamma distribution. However, this is difficult and may not be possible using solely a hypoexponential distribution, as is known from previous literature (Section~\ref{Sec:3moment}). Effectively, it may appear that the approximation is easily generalized to match more than two moments. However, we now show that the roots of the corresponding polynominal can not be real-valued and thus prove
\begin{theorem}\label{Thm:SmoothedThreeMomentAppendix}
Let $\mathcal{X}$ be a $\GammaDist(j,a)$-distributed random variable where $j \notin \mathbb{N}.$ Then, it is not possible to match three or more moments of $\mathcal{X}$ using the smoothed hypoexponential distribution.
\end{theorem}

Before formulating and proving this result, we first mention two facts about the cumulant-generating functions of Gamma and hypoexponential random variables. The cumulant generating function of a gamma distributed random variable $X$ with shape and rate parameters $j,a$ is given by
\begin{equation*}
K_X(\theta) = -j \log(1-\theta/a) = j\sum_{m=1}^{\infty} \frac{(\theta/a)^m}{m}  
\end{equation*}
Therefore, the cumulants $\kappa^X_m = \left.\frac{\d^m}{\d\theta^m} K_X(\theta) \right|_{\theta=0}$ are equal to
\begin{equation*}
 \kappa^X_m = \frac{j}{a^m} (m-1)!.
\end{equation*}
Let $Y$ be a random variable with a hypoexponential distribution
with rate parameters $a_1,\dots,a_n$. The cumulant generating function of $Y$ is given by
\begin{equation*}
 K_Y(\theta) = -\sum_{k=1}^n \log(1-\theta/a_k) = \sum_{m=1}^{\infty} \frac{\theta^m}{m}  \sum_{k=1}^n \frac{1}{a_k^m}
\end{equation*}
Therefore, the cumulants of $Y$ are given by
\begin{equation*}
 \kappa_m^Y = (m-1)! \sum_{k=1}^n \frac{1}{a_k^m}.
\end{equation*}

Now, suppose that we want to match $m$ moments of $X$ and $Y_s$ by choosing $m <n$ rates as free variables. Of the $n$ rates in the smoothed hypoexponential approximation, the final $n-m$ rates $\lambda_k$ will be equal to $a= j/\tau$. Set $x_k = 1/\lambda_k$ and note that, without loss of generality, we can assume that $a = 1$, because otherwise we can replace $X$ and $Y$ by $aX$ and $aY_s$, respectively. Hence, the cumulants of $X$ are given by $\kappa_m^X = j (m-1)!$. 

The unknown $m$ rates satisfy the following system of equations.
\begin{equation*}
    \left((n-m) \cdot 1 + x_1^k + \cdots + x_m^k\right) (k-1)! = \kappa_k^Y = \kappa_k^X =  j(k-1)!\,,\quad k = 1,\dots,m
\end{equation*}
which can be written as
\begin{equation}\label{eqn:m-duration-system}
 x_1^k + x_2^k + \cdots + x_m^k = (m-1+\{j\})\,,\quad k = 1,\dots,m
\end{equation}
where $\{j\}$ is the fractional part of $j$. The strategy is to construct a polynomial $f_m$ of which the $x_i$ are the roots. By construction this polynomial equals
\begin{equation*}
 f_m(x) = \prod_{k=1}^m (x-x_k) = \sum_{k=0}^m (-1)^ke_{k}(x_1,\dots, x_m) x^{m-k}
\end{equation*}
where the $e_k$ are the elementary symmetric polynomials.

\begin{theorem}\label{Thm:fm}
	Write $(z)_k = z (z-1)(z-2)\cdots (z-k+1)$ for the falling Pochhammer symbol. 
	The polynomial $f_m(x)$ is equal to
	\begin{equation*}
	 f_m(x) = \sum_{k=0}^m (-1)^k x^{m-k} \frac{(m-1+\{j\})_k}{k!}
	\end{equation*}
\end{theorem}
\begin{proof}
	Write $p_k(x_1,\dots,x_m) = \sum_{i=1}^m x_i^k$ for the power sum of the roots $x_i$.
	According to Newton's identity for symmetric polynomials, we have
	\begin{equation*}
	 k e_k(x_1,\dots, x_m) = \sum_{\ell=1}^k (-1)^{\ell-1} e_{k-\ell}(x_1,\dots,x_m) p_{\ell}(x_1,\dots, x_m)
	\end{equation*}
	According to Eq~\ref{eqn:m-duration-system}, we know that the power sums equal
	\begin{equation*}
	p_{k}(x_1,\dots, x_m) = (m-1 + \{j\})
	\end{equation*}	
	Hence, in order to prove Theorem~\ref{Thm:SmoothedThreeMomentAppendix}, it remains to show that the coefficients of $f_m$ satisfy Newtons identity, i.e. we have to show that 
	\begin{equation}\label{eqn:hypothesis}
	k \frac{(m-1+\{j\})_k}{k!} = \sum_{\ell = 1}^k (-1)^{\ell-1} \frac{(m-1+\{j\})_{k-\ell}}{(k-\ell)!} (m-1+\{j\}) 
	\end{equation}
	According to the Chu-Vandermonde identity, we have
	\begin{equation*}
	 (a+b)_{N} = \sum_{K=0}^N \binom{N}{K} (a)_{K}(b)_{N-K}
	\end{equation*}
	By taking $a = -1$, $b = m-1+\{j\}$, and $N = k-1$, we find that 
	\begin{equation*}
	 \sum_{\ell=1}^{k} \binom{k-1}{\ell-1} (m-1+\{j\})_{k-\ell} (-1)_{\ell-1} = (m-2+\{j\})_{k-1}
	\end{equation*}
	which we can re-write as
	\begin{equation*} 
	(m-1+\{j\})_{k} = (m-1+\{j\})\sum_{\ell=1}^k \frac{(-1)_{\ell-1}(k-1)!}{(k-\ell)! (\ell-1)!} (m-1+\{j\})_{k-\ell}
	\end{equation*}
	From which we see that Eq~\ref{eqn:hypothesis} is indeed true.
\end{proof}

We now introduce a related family of polynomials $\{g_m(x) \}_{m=0}^{\infty}$ defined by 
\begin{align*}
    g_m(x) = x^m f_m(1/x)
\end{align*}
and investigate the real roots of the $g_m$.

\begin{lemma}\label{lem:gm-deriv}
	Let $g_m(x) = x^m f_m(1/x)$, then $g_0(x) \equiv 1$ and for $m > 0$ we have
	\begin{equation*}
	 \frac{\d}{\d x}g_{m}(x) = - (m -1 + \{j\})g_{m-1}(x)
	\end{equation*}
	In addition, for all $m \geq 0$, we have $g_m(0) = 1$.
\end{lemma}

\begin{proof}
    By definition, $g_m(0) = 1$ for all $m$ and $g_0(x) \equiv 1$.
	By differentiating $g_m$, we easily find that 
	\begin{equation*}
	\begin{split}
	 \frac{\d}{\d x} g_m(x) &= \sum_{k=0}^m (-1)^k\frac{(m-1+\{j\})_{k}}{k!} k x^{k-1} \\
	 &= -(m-1+\{j\})\sum_{k=1}^m (-1)^{k-1}\frac{(m-1+\{j\})_{k-1}}{(k-1)!} x^{k-1} \\
	 &= -(m-1+\{j\}) g_{m-1}(x)
	 \end{split}
	\end{equation*}
	which proves the Lemma.
\end{proof}

\begin{lemma}\label{lem:gm1}
	for all $m>0$, we have $g_m(1) = (-1)^m\binom{m-2+\{j\}}{m}$.
	In particular $g_m(1) < 0$ for even $m$ and $g_m(1) > 0$ for odd $m$.
\end{lemma}

\begin{proof}
	This follows directly from a well known fact about alternating sums of binomial coefficients, which we will prove for the reader's convenience. We first notice that Pascal's rule generalizes to binomial coefficients with a non-integer top argument:
	\begin{equation}\label{eqn:pascal}
	 \binom{\alpha}{m} + \binom{\alpha}{m+1} = \binom{\alpha+1}{m+1}, \quad \forall \alpha \in \mathbb{C}, m \in \mathbb{Z}_{\geq 0}
	\end{equation}
	We have to prove that 
	\begin{equation}\label{eqn:alt-binom}
	 \sum_{k=0}^m(-1)^k \binom{\alpha}{k} = (-1)^m\binom{\alpha-1}{m}
	\end{equation}
	For $m = 0$, this is trivially true.
	Suppose that Eq~\ref{eqn:alt-binom} holds for all $\ell < m$, then
	\begin{multline*}
	\sum_{k=0}^m(-1)^k \binom{\alpha}{k} 
	= (-1)^{m-1} \binom{\alpha-1}{m-1} + (-1)^{m} \binom{\alpha}{m} \\
	= (-1)^{m} \left[ \binom{\alpha-1}{m-1} + \binom{\alpha-1}{m} - \binom{\alpha-1}{m-1} \right]
	= (-1)^m \binom{\alpha-1}{m}
	\end{multline*}
	Now choosing $\alpha = m-1+\{j\}$ completes the lemma.
\end{proof}

\begin{lemma}\label{lem:gm-lb}
If  $m$ is odd and $0 < x < 1$, then $g_m(x) > (1-x)^m$.
\end{lemma}
\begin{proof}
	As $(1-x)^{m-1+\{j\}} = \sum_{k=0}^{\infty} (-1)^k x^k \frac{(m-1+\{j\})_k}{k!}$, we find that
	\begin{equation*}
	 (1-x)^{m-1+\{j\}} - g_m(x) = \sum_{k=m+1}^{\infty} (-1)^k x^k \frac{(m-1+\{j\})_k}{k!}
	\end{equation*}
	We will show that all coefficients of the power series on the right-hand side are negative.
	Let $k = m+\ell$ with $\ell > 0$, then as $m$ is odd, we have
	\begin{equation}
	 (-1)^k (m-1+\{j\})_k = -(-1)^{\ell} (m-1+\{j\})_m (\{j\}-1)_{\ell}
	 \label{Eq:PochhammerSigns}
	\end{equation}
	The first Pochhammer symbol on the right-hand side of \eqref{Eq:PochhammerSigns} is positive, and the second has the same sign as $(-1)^{\ell}$. This proves the lemma.
\end{proof}

\begin{lemma}\label{lem:gm-roots-01}
	Let $j \not\in \mathbb{Z}$ and $m > 0$. Then, $g_m(x)$  has zero roots in the unit interval $[0,1]$ if $m$ is odd and $g_m(x)$ has one root in the unit interval $[0,1]$ if $m$ is even.  
\end{lemma}
\begin{proof}
	We know from lemma~\ref{lem:gm-lb} that for odd $m$, we have $
	g_m(x) > (1-x)^{m-1+\{j\}} > 0$ thus showing the first conclusion of the claim. 
	
	Now, suppose that $m $ is even. Then $\tfrac{\d}{\d x} g_m(x) = -(m-1+\{j\})g_{m-1}(x) < 0$,
	because $m-1$ is odd. Hence, $g_m(x)$ is monotone decreasing. Further, $g_m(0) = 1$ and $g_m(1) < 0$ ( by lemma~\ref{lem:gm1}), so $g_m(x)$ must have exactly one root on $[0,1]$.
\end{proof}

\begin{lemma}\label{lem:gm-roots-1infty} 
 Assume that $j \not\in \mathbb{Z}$ and $m > 0$, then $g_m(x)$ has exactly $1$ root in $(1,\infty)$.
\end{lemma}
\begin{proof}
We proceed by induction on $m$.  For $m = 1$, the polynomial $g_m(x) = 1- \{j\}x$ is linear and has exactly $1$ root equal to  $1/\{j\} > 1$. 

Now, suppose that the lemma is true for $\ell < m$, and first consider the case where $m$ is even. From lemma~\ref{lem:gm1}, we know that $g_m(1) < 0$. As the coefficient of $x^m$ in $g_m(x)$ is positive, we also know that $g_m(x) > 0$ for large enough $x$. Hence, $g_m(x)$ has at least one root in $(1,\infty)$. 
 
 Now suppose that $g_m(x)$ has more than one root on $(1,\infty)$. As again $g_m(x) > 0$ for large enough $x$, the number of roots must be odd and $\geq 3$. Consequently, the number of local extrema of $g_m(x)$ for $x \in (0,\infty)$ must be $\geq 2$. As $\tfrac{\d}{\d x}g_m(x) = -(m-1+\{j\})g_{m-1}(x)$ (lemma~\ref{lem:gm-deriv}), this would mean that $g_{m-1}(x)$ has more than $1$ root on $(1,\infty)$, which is impossible according to the induction hypothesis. Hence, $g_m(x)$ has exactly one root on $(1,\infty)$.
 
 The $m$ odd case follows from an identical argument applied to $-g_m(x)$ and establishes the claim.
\end{proof}

\begin{theorem}
	Assume that $j \not\in \mathbb{Z}$ and $m>0$. Then, the polynomial $f_m(x)$ has at least one and at most two roots in $\mathbb{R}$.
\end{theorem}
\begin{proof}
	From lemma~\ref{lem:gm-roots-01} and lemma~\ref{lem:gm-roots-1infty} it follows that 
	$f_m$ has one or two real, positive roots. From Theorem~\ref{Thm:fm}, we also see that $f_m$ has no non-positive roots. This proves the Theorem.
\end{proof}

\end{document}